# A Synthetic Version of Lie's Second Theorem

Matthew Burke[*]

May 24, 2016

## Abstract

We formulate and prove a twofold generalisation of Lie's second theorem that integrates homomorphisms between formal group laws to homomorphisms between Lie groups. Firstly we generalise classical Lie theory by replacing groups with categories. Secondly we include categories whose underlying spaces are not smooth manifolds. The main intended application is when we replace the category of smooth manifolds with a well-adapted model of synthetic differential geometry. In addition we provide an axiomatic system that provides all the abstract structure that is required to prove Lie's second theorem. As a part of this abstract structure we define the notion of enriched mono-coreflective subcategory which makes precise the notion of a subcategory of local models.

## Contents

[*]The author acknowledges the support of an International Macquarie University Research Excellence Scholarship.

# 1 Introduction

In classical Lie theory there is an adjunction

$$FGLaw \underset{(-)_\infty}{\overset{(-)_{int}}{\rightleftarrows}} LieGrp \qquad (1)$$

between the category $FGLaw$ of formal group laws and the category $LieGrp$ of Lie groups. When we restrict the domain of $(-)_\infty$ to the simply connected Lie groups, Lie's second theorem tells us that $(-)_\infty$ is full and faithful and Lie's third theorem tells us that $(-)_\infty$ is essentially surjective. We refer to [17] for the classical theory of Lie groups, Lie algebras and formal group laws. For instance we can combine Theorem 3 of Section V.6 and Theorem 2 of Section V.8 of Part 2 in [17] to obtain the equivalence above. Given a Lie group $\mathbb{G}$ we think of the formal group law $\mathbb{G}_\infty$ as consisting of all the data contained in an infinitesimal neighbourhood of the identity element of $\mathbb{G}$ and so the functor $(-)_{int}$ is interpreted as specifying a way to extend local data to global data.

We now recall some definitions from the established multi-object generalisation of Lie theory involving Lie groupoids. A *Lie groupoid* is a groupoid in the category $Man$ of smooth manifolds such that the source and target maps are submersions. A *Lie algebroid* is a vector bundle $A \to M$ together with a bundle homomorphism $\rho : A \to TM$ such that the space of sections $\Gamma(A)$ is a Lie algebra satisfying the following Leibniz law: for all $X, Y \in \Gamma(A)$ and $f \in C^\infty(M)$ the equality

$$[X, fY] = \rho(X)(f) \cdot Y + f[X, Y]$$

holds. In multi-object Lie theory we have a functor

$$LieAlgd \xleftarrow{(-)_\infty} LieGpd$$

from the category of Lie groupoids to the category of Lie algebroids which is full and faithful but not essentially surjective when we impose the appropriate connectedness conditions. Any Lie algebroid integrates to a topological groupoid, its Weinstein groupoid [5], but there can be obstructions to putting a smooth structure on it.

In this paper we generalise this situation in two main ways. Firstly we replace groupoids with categories. Secondly we replace the category $Man$ with a well-adapted model $\mathcal{E}$ of synthetic differential geometry as described in Section 2.1. The category $\mathcal{E}$ contains the category of smooth manifolds (with boundary) as a full subcategory but unlike $Man$ the category $\mathcal{E}$ is closed under all limits and colimits and contains rigorously defined infinitesimal objects. This allows us to replace the Lie correspondence (1) with a correspondence between two types of categories in $\mathcal{E}$. The first type of category (called jet categories) are those for which every arrow is infinitesimally close to an identity arrow. The jet categories play an analogous role to the role played by the formal group laws in classical Lie theory. The second type of category (called integral complete categories) are those in which every time-dependent left-invariant vector field admits a local solution. The integral complete categories play an analogous role to the role played by the Lie groups in the classical theory.

In Section 4.2 and Section 6.2 we show that the subcategories $Cat_\infty(\mathcal{E})$ of jet categories and $Cat_{int}(\mathcal{E})$ of integral complete categories are coreflective and reflective respectively in $Cat(\mathcal{E})$. The adjunction that we use to replace the Lie adjunction (1) is the composite of the reflection and the coreflection:

$$Cat_\infty(\mathcal{E}) \xrightarrow[\underset{(-)_\infty}{\perp}]{i} Cat(\mathcal{E}) \xrightarrow[\underset{j}{\perp}]{(-)_{int}} Cat_{int}(\mathcal{E}) \qquad (2)$$

The main result of this paper is that if we impose certain enriched connectedness conditions (which are described in Section 7) then the functor $(-)_\infty j$ is full and faithful. In fact the construction the right hand side of (2) from the left hand side at a higher level of generality. In Definition 4.9 we give the definition of $\mathcal{E}$-mono-coreflective subcategory which picks out the key structural components that (along with a choice of interval object) allow us to prove Lie's second theorem. Note that in this case $\mathcal{E}$ can be an arbitrary topos. We think of an $\mathcal{E}$-mono-coreflective subcategory as a subcategory of local models; fixing an $\mathcal{E}$-mono-coreflective subcategory amounts to defining the term 'local model'. This means that we can construct a version of Lie theory and prove Lie's second theorem in any situation where we can construct an $\mathcal{E}$-mono-coreflective subcategory.

This new correspondence is useful in several ways. Firstly it applies to groupoids whose underlying space is non-classical in nature. This avoids the problems existing in the literature (see for instance [5]) concerning the non-integrability of Lie algebroids. Secondly it provides a novel analytic approximation to a Lie groupoid. Invariably Lie groupoids have been approximated by Lie algebroids which constitute linear approximations. As such it would be interesting to compare the objects of $Cat_\infty(\mathcal{E})$ to the formal symplectic groupoids in [4].

Finally the categorical nature of the constructions and results in this paper make them amenable to generalisation and application in other areas. For instance the constructions can be carried out without additional difficulties in the big Zariski topos (see Section 4 in Chapter 3 of [14]) and hence provides a candidate for a formulation of Lie theory that involves the (not necessarily smooth) group schemes defined therein. The advantage of working relative to an $\mathcal{E}$-mono-coreflective subcategory is that we can separate the local approximation procedure from the global and higher dimensional aspects of the theory. This means that we can carry out the local approximation of groupoids and categories separately in Section 4.1 and Section 4.3 respectively but then use the same global constructions to complete the proof of Lie's second theorem. Furthermore in Section 7.2 we explicitly arrange the higher dimensional data required to prove Lie's second theorem into truncated cubical objects. This makes the formulation of Lie theory in this paper a promising first step in establishing a Lie theory for higher categories and groupoids in the area of derived algebraic geometry and $(\infty, 1)$-categories.

**Relationship to Classical and Multi-object Lie Theory**

In [2] we justify the general constructions in this paper by explaining how they relate to the classical ones in the case that $\mathcal{E}$ is a well-adapted model of synthetic differential geometry. For instance in [2] we show that when we restrict to the subcategory of $Cat(\mathcal{E})$ consisting of the classical Lie groups the functor $(-)_\infty$ corresponds to the formal group law construction given in the Introduction to [7]. In multi-object Lie theory we restrict attention to the groupoids that are source simply connected (i.e. have simply connected source fibres). For arbitrary groupoids in $\mathcal{E}$ the enriched path and simply connectedness conditions that we describe in Section 7 are slightly stronger than their classical counterparts. However in [2] we show that a Lie groupoid is path connected iff it is enriched path connected and that a Lie groupoid is simply connected iff it is enriched simply connected.

In Section 4.3 we show how to specialise our results so that they are phrased in terms of groupoids rather than categories. To do this we need to assume that the neighbour relation defined in Section 3.2 is symmetric on the arrow space of our groupoids. Therefore in [2] we justify this by showing that the infinitesimal neighbour relation is symmetric for any Lie groupoid. In Section 7.5 we see that in order to prove Lie's second theorem in this context we need an additional assumption that is invisible in the classical theory: that the jet part $\mathbb{C}_\infty$ of our category is enriched path connected. In [2] we show that for any Lie groupoid $\mathbb{G}$ the groupoid $\mathbb{G}_\infty$ is enriched path connected.

## 2 Preliminaries

In this section we give a brief overview of the parts of the theory of synthetic differential geometry and enriched factorisation systems which will be useful to us. In Section 2.1 describe the set $Spec(Weil)$ of nilpotent infinitesimals which we will use to construct our infinitesimal approximations. In Section 2.2 we recall two different ways to generate an enriched factorisation system. The first uses the method of wide intersections to generate a factorisation system from arrows in the right class; the second uses the small object argument to generate a factorisation system from arrows in the left class.

### 2.1 Synthetic Differential Geometry

In synthetic differential geometry we replace the category $Man$ of smooth manifolds (with boundary) with a certain kind of Grothendieck topos $\mathcal{E}$ called a well-adapted model of synthetic differential geometry. We now describe a few of the key properties of $\mathcal{E}$. Firstly there is a full and faithful embedding $\iota : Man \rightarrowtail \mathcal{E}$ and therefore a ring $R = \iota\mathbb{R}$ and unit interval $I = \iota I$ in $\mathcal{E}$. In addition we have the objects

$$D_k = \{x \in R : x^{k+1} = 0\}$$

which are not terminal. In fact the fundamental Kock-Lawvere axiom holds: the arrow $\alpha : R^{k+1} \to R^{D_k}$ defined by

$$(a_0, a_1, ..., a_k) \mapsto (d \mapsto a_0 + a_1 d + ... + a_k d^k)$$

is an isomorphism. We write $D_\infty = \bigcup_i D_i$ and $D = D_1$. Using the Kock-Lawvere axiom we can show that $\iota(TM) \cong (\iota M)^D$ as vector bundles over $\iota M$ and that the Lie bracket corresponds to an infinitesimal commutator. For more details on this construction and synthetic differential geometry in general see [12]. A set of non-classical objects that will be useful in the sequel are the Weil spectra which are of the form:

$$Spec(Weil) = \left\{ (x_1, ..., x_n) : \bigwedge_{i=1}^{n}(x_i^{k_i} = 0) \wedge \bigwedge_{j=1}^{m}(p_j = 0) \right\}$$

where $n, m \in \mathbb{N}_{\geq 0}$, $k_i \in \mathbb{N}_{>0}$ and the $p_j$ are polynomials in the $x_i$.

#### 2.1.1 The Amazing Right Adjoint

An important property of spectra of Weil algebras is that they are 'atomic' objects of the topos. In short this says that they are small enough to only fit in one summand of any structure that we construct by glueing together other smaller structures. The next definition makes this idea precise.

**Definition 2.1.** An object $X$ in a category $\mathcal{E}$ is *atomic* iff the endofunctor

$$\mathcal{E} \xrightarrow{(-)^X} \mathcal{E}$$

defined using the internal hom has a right adjoint.

**Proposition 2.2.** *The object $D$ is atomic for all $D \in Spec(Weil)$ in a well-adapted model of synthetic differential geometry.*

*Proof.* This follows from the Example in Appendix 4 of [15]. $\square$

## 2.2 Enriched Factorisation Systems

In this section we sketch the theory of enriched factorisation systems that we use to construct the jet part of a category. We follow [6] by defining the orthogonality of arrows in terms of hom-objects. Although we work analogously to the treatment of weak enriched factorisation systems in [16] we mainly make use of the account of (orthogonal) enriched factorisation systems that is [13]. We refer to [11] for the basic concepts of enriched category theory.

**Notation 2.3.** Let $\mathcal{V}$ be a monoidal category. Let $\mathcal{C}$ be a $\mathcal{V}$-category. Then we write $\mathcal{C}_0$ for the underlying ordinary category of $\mathcal{C}$.

**Definition 2.4.** The arrow $l$ is left $\mathcal{V}$-orthogonal to $r$ (written $l \perp_\mathcal{V} r$) iff

$$\begin{array}{ccc} \mathcal{C}(B, X) & \xrightarrow{\mathcal{C}(l,X)} & \mathcal{C}(A, X) \\ {\scriptstyle \mathcal{C}(B,r)} \downarrow & & \downarrow {\scriptstyle \mathcal{C}(A,r)} \\ \mathcal{C}(B, Y) & \xrightarrow{\mathcal{C}(l,Y)} & \mathcal{C}(A, Y) \end{array}$$

is a pullback in $\mathcal{V}$.

**Definition 2.5.** Let $S$ be a class of arrows in $\mathcal{C}_0$. Then the right $\mathcal{V}$-orthogonal complement of $S$ is the class:

$$S^{\perp_\mathcal{V}} := \{f \in \mathcal{C}_0^{\mathbf{2}} : \forall s \in S.\ s \perp_\mathcal{V} f\}$$

and the left $\mathcal{V}$-orthogonal complement of $S$ is the class:

$$^{\perp_\mathcal{V}}S := \{f \in \mathcal{C}_0^{\mathbf{2}} : \forall s \in S.\ f \perp_\mathcal{V} s\}$$

**Definition 2.6.** The pair $(L, R)$ is a $\mathcal{V}$-prefactorisation system on $\mathcal{C}$ iff $L^{\perp_\mathcal{V}} = R$ and $L = {}^{\perp_\mathcal{V}}R$.

**Definition 2.7.** The pair $(L, R)$ is a $\mathcal{V}$-factorisation system on $\mathcal{C}$ iff $(L, R)$ is a prefactorisation system and $(L, R)$-factorisations exist: i.e. for every $f \in \mathcal{C}_0^{\mathbf{2}}$ there exist $l \in L, r \in R$ such that $f = r \circ l$.

The next result provides sufficient conditions for a pair $(L, R)$ to be a $\mathcal{V}$-factorisation system. We use these conditions in Section 3.2 when defining the jet factorisation system using an infinitesimal neighbour relation.

**Lemma 2.8.** *The pair $(L, R)$ is a $\mathcal{V}$-factorisation system iff*

1. *the classes $L$ and $R$ are replete,*
2. *if $l \in L$ and $r \in R$ then $l \perp_{\mathcal{V}} r$,*
3. *for every map $f$ in $\mathcal{C}$, there exist $f_r \in R$ and $f_l \in L$ such that $f = f_r f_l$.*

Recall the following result for generating an enriched factorisation system from its right class. We use this result to define the jet factorisation system in Section 3.1. It is Lemma 3.1 of [3] where a sketch of the proof is given and Proposition 7.1 in [13] where a full proof is given.

**Proposition 2.9.** *Let $R$ be a class of arrows in a category $\mathcal{C}$. Suppose that $R$ is contained in the class of monomorphisms, is closed under composition and contains all the isomorphisms. Suppose that the pullback of an arrow in $R$ along an arbitrary arrow in $\mathcal{C}$ exists in $\mathcal{C}$ and is again in $R$. Suppose further that all intersections of arrows in $R$ exist in $\mathcal{C}$ and are again in $R$. Then $(^{\perp}R, R)$ is a factorisation system on $\mathcal{C}$.*

Now we recall a way to generate an enriched factorisation system from a generating set of arrows in the left class. We use this result to define the integral factorisation system in Section 6.1.

**Proposition 2.10.** *Let $\mathcal{C}$ be a $\mathcal{V}$-category such that its underlying category $\mathcal{C}_0$ is locally presentable. Let $\Sigma$ be a set of arrows in $\mathcal{C}$. Then there is a factorisation system $(L, R)$ on $\mathcal{C}_0$ such that $R = \Sigma^{\perp_{\mathcal{V}}}$ and $L = {}^{\perp_{\mathcal{V}}}R$.*

# 3 The Jet Factorisation System

We construct the infinitesimal (or jet) part of a category using an enriched factorisation system that we call the jet factorisation system. We give two characterisations of the jet factorisation system. The first in Section 3.1 is more direct and is used to prove that the jet part of a category is closed under composition. The second is defined in terms of an infinitesimal neighbour relation and is used when we need a more concrete description of the left and right classes. For instance we use the neighbour relation in Section 3.3 to find conditions under which the left class is stable under pullback and in Section 4.3 to find a necessary and sufficient condition for the jet part of a groupoid to be a groupoid (rather than just a category).

## 3.1 Jet Factorisation in the Slice Topos

We define the jet factorisation system on any slice category $\mathcal{E}/M$ of the well-adapted model of synthetic differential geometry $\mathcal{E}$. Since it is a topos the category $\mathcal{E}$ is locally cartesian closed. Furthermore we show that for any arrow $f : X \to Y$ in $\mathcal{E}$ both the pullback functor $f^* : \mathcal{E}/Y \to \mathcal{E}/X$ and its left adjoint $\Sigma_f : \mathcal{E}/X \to \mathcal{E}/Y$ preserve the left class of the jet factorisation systems on $\mathcal{E}/Y$ and $\mathcal{E}/X$ respectively. This will be used in the next section to define the composition operation on the jet part of a category in $\mathcal{E}$. In the case $M = 1$ the right class of the jet factorisation system has been studied before. It is the class of formal-etale maps in I.17 of [12] and the class of formally-open morphisms in Section 1.2 of Volume 3 of [8]. For the standard theory of toposes we refer to [14].

In this section $\mathcal{E}$ will be a smooth topos and $M$ an object of $\mathcal{E}$. To begin with let us recall the definition of slice category. It can be found for example in construction 4 of Section 1.6 in [1].

**Definition 3.1.** The *slice category* $\mathcal{E}/M$ *of a category* $\mathcal{E}$ *over an object* $M \in \mathcal{E}$ has as objects all arrows $f \in \mathcal{E}$ such that the codomain of $f$ is $M$. To keep track of the domain of $f$ we write the objects of $\mathcal{E}/M$ in the form $(dom(f), f)$. An arrow $g : (X, f) \to (X', f')$ in $\mathcal{E}/M$ is an arrow $g : X \to X'$ in $\mathcal{E}$ such that $f' \circ g = f$.

The following is part of Theorem 1.42 in [9].

**Theorem 3.2.** *Let $\mathcal{E}$ be a topos, $X$ an object of $\mathcal{E}$. Then $\mathcal{E}/X$ is a topos.*

Now we define the jet factorisation system.

**Definition 3.3.** An arrow $r : X \to Y$ in $\mathcal{E}/M$ is *jet closed* iff it is a monomorphism and

$$\begin{array}{ccc} X^{(M \times D, \pi_1)} & \xrightarrow{X^{(1_M, 0)}} & X \\ \downarrow{r^{(M \times D, \pi_1)}} & & \downarrow{r^{(M, 1_M)}} \\ Y^{(M \times D, \pi_1)} & \xrightarrow{Y^{(1_M, 0)}} & Y \end{array}$$

is a pullback in $\mathcal{E}/M$ for all $D$ in $Spec(Weil)$. An arrow $l : A \to B$ in $\mathcal{E}/M$ is *jet dense* iff for all jet closed arrows $r$ the square

$$\begin{array}{ccc} X^B & \xrightarrow{X^l} & X^A \\ \downarrow{r^B} & & \downarrow{r^A} \\ Y^B & \xrightarrow{Y^l} & Y^A \end{array}$$

is a pullback in $\mathcal{E}/M$.

**Definition 3.4.** The *jet factorisation system on* $\mathcal{E}/M$ is the pair $(L_\infty, R_\infty)$ where $L_\infty$ is the class of jet dense arrows and $R_\infty$ is the class of jet closed arrows.

**Remark 3.5.** The fact that $(L_\infty, R_\infty)$ is an $\mathcal{E}/M$ factorisation system follows from Proposition 2.9.

We now relate the jet factorisation systems on different slice categories of $\mathcal{E}$ by using the fact that $\mathcal{E}$ is locally cartesian closed.

**Proposition 3.6.** *Let $f : G \to M$ be an arrow in $\mathcal{E}$. Let $f^* : \mathcal{E}/M \to \mathcal{E}/G$ be the functor defined by pullback along $f$. Then $f^*$ preserves exponentials and has both a left adjoint $\Sigma_f$ and right adjoint $\Pi_f$; the left adjoint $\Sigma_f$ is given by postcomposition with $f$.*

$$\mathcal{E}/G \underset{\Pi_f}{\overset{\Sigma_f}{\underset{\longrightarrow}{\overset{\longrightarrow}{\leftarrow f^* -}}}} \mathcal{E}/M$$

*Proof.* This is Theorem 2 on page 193 in [14]. $\square$

**Lemma 3.7.** *Let $\rho : X \rightarrowtail Y$ be a jet closed arrow in $\mathcal{E}/M$ and $f : G \to M$ an arrow in $\mathcal{E}$. Then $f^*(\rho)$ is a jet-closed arrow in $\mathcal{E}/G$.*

*Proof.* Since $\rho$ is jet closed in $\mathcal{E}/M$ we have that for all $D \in Spec(Weil)$ the following square is a pullback:

$$\begin{array}{ccc} X^{(M \times D, \pi_1)} & \xrightarrow{X^{(1_M, 0)}} & X \\ \downarrow{\rho^{(M \times D, \pi_1)}} & & \downarrow{\rho} \\ Y^{(M \times D, \pi_1)} & \xrightarrow{Y^{(1_M, 0)}} & Y \end{array}$$

Using the fact that $f^*$ preserves exponentials we see that:

$$f^* \left( \begin{array}{ccc} X^{(M \times D, \pi_1)} & \xrightarrow{X^{(1_M, 0)}} & X \\ \downarrow{\rho^{(M \times D, \pi_1)}} & & \downarrow{\rho} \\ Y^{(M \times D, \pi_1)} & \xrightarrow{Y^{(1_M, 0)}} & Y \end{array} \right) \cong \begin{array}{ccc} f^*(X)^{(G \times D, \pi_1)} & \xrightarrow{f^*(X)^{(1_G, 0)}} & f^*(X) \\ \downarrow{f^*(\rho)^{(G \times D, \pi_1)}} & & \downarrow{f^*(\rho)} \\ f^*(Y)^{(G \times D, \pi_1)} & \xrightarrow{f^*(Y)^{(1_G, 0)}} & f^*(Y) \end{array}$$

Then using the fact that $f^*$ is a right adjoint we deduce that the right hand square is a pullback for all $D \in Spec(Weil)$ and so $f^*(\rho)$ is jet-closed in $\mathcal{E}/G$. $\square$

**Lemma 3.8.** *Let $F \dashv U$ be adjoint functors. Suppose that $F$ preserves products. Then:*

$$(UA)^B \cong U(A^{FB})$$

*Proof.* We will establish a natural bijection between the generalised elements of both sides:

$$\frac{\dfrac{\dfrac{\dfrac{X \to (UA)^B}{X \times B \to UA}}{F(X \times B) \to A}}{FX \to A^{FB}}}{X \to U(A^{FB})}$$

as required. $\square$

**Lemma 3.9.** *Let $\rho : X \rightarrowtail Y$ be a jet closed arrow in $\mathcal{E}/G$ and $f : G \to M$ an arrow in $\mathcal{E}$. Then $\Pi_f(\rho)$ is a jet-closed arrow in $\mathcal{E}/M$.*

*Proof.* Since $\rho$ is jet closed in $\mathcal{E}/G$ we have that for all $D \in Spec(Weil)$ the following square is a pullback:

$$\begin{array}{ccc} X^{(G \times D, \pi_1)} & \xrightarrow{X^{(1_G, 0)}} & X \\ \downarrow{\rho^{(G \times D, \pi_1)}} & & \downarrow{\rho} \\ Y^{(G \times D, \pi_1)} & \xrightarrow{Y^{(1_G, 0)}} & Y \end{array}$$

Using Lemma 3.8 we see that:

$$\Pi_f \left( \begin{array}{ccc} X^{(G \times D, \pi_1)} & \xrightarrow{X^{(1_G, 0)}} & X \\ \downarrow{\rho^{(G \times D, \pi_1)}} & & \downarrow{\rho} \\ Y^{(G \times D, \pi_1)} & \xrightarrow{Y^{(1_G, 0)}} & Y \end{array} \right) \cong \begin{array}{ccc} \Pi_f(X)^{(M \times D, \pi_1)} & \xrightarrow{\Pi_f(X)^{(1_M, 0)}} & \Pi_f(X) \\ \downarrow{\Pi_f(\rho)^{(M \times D, \pi_1)}} & & \downarrow{\Pi_f(\rho)} \\ \Pi_f(Y)^{(M \times D, \pi_1)} & \xrightarrow{\Pi_f(Y)^{(1_M, 0)}} & \Pi_f(Y) \end{array}$$

Then using the fact that $\Pi_f$ is a right adjoint we deduce that the right hand square is a pullback for all $D \in Spec(Weil)$ and so $\Pi_f(\rho)$ is jet-closed in $\mathcal{E}/M$. $\square$

**Corollary 3.10.** *Let $l$ be jet dense in $\mathcal{E}/G$ and $f : G \to M$ an arrow in $\mathcal{E}$. Then $\Sigma_f(l)$ is jet dense in $\mathcal{E}/M$.*

**Corollary 3.11.** *Let $\lambda$ be a jet dense arrow in $\mathcal{E}/M$ and $f : G \to M$ an arrow in $\mathcal{E}$. Then $f^*(\lambda)$ is jet dense in $\mathcal{E}/G$.*

### 3.2 Jet Factorisation Using Neighbours

The jet factorisation system presented in Section 3.1 can be thought of as a 'perturbation' of the standard $(Epi, Mono)$-factorisation in a topos. Intuitively speaking, if $f : A \to B$ is a jet dense arrow and $b$ is an element of $B$ then although there might not exist an element $a$ of $A$ such that $fa = b$ there does exist an element $a'$ of $A$ such that $fa'$ is 'infinitesimally close' to $b$. We can give a similar heuristic description for the jet closed arrows. If $g : X \rightarrowtail Y$ is a jet closed arrow then it is a monomorphism by definition. But $g$ satisfies an additional condition: if $x$ is an element of $X$ and $y$ is an element of $Y$ such that $gx$ is infinitesimally close to $y$ then there exists an element $x'$ in $X$ such that $gx' = y$. In this section we make these ideas precise by defining a reflexive relation $\sim$ in the internal logic of the topos $\mathcal{E}/M$ for which $a \sim b$ encodes the idea that $b$ is contained in some infinitesimal perturbation (or jet) which is based at $a$. Then we define a factorisation system using this relation which corresponds to our intuitive idea of perturbing the $(Epi, Mono)$-factorisation in $\mathcal{E}/M$. Finally we show that this factorisation system in fact coincides with the jet factorisation system.

First we recall the definition of generalised element in a category from Definition 1.1 in Part II of [12].

**Definition 3.12.** Let $R$ be an object in a category $\mathcal{E}$. A *generalised element* of $R$ is an arrow in $\mathcal{E}$ with codomain $R$. The domain of the arrow is called the *stage of definition* of the element.

**Notation 3.13.** We write $r \in_X R$ to denote that $r$ is an arrow $X \to R$ in $\mathcal{E}$ and hence $r$ is an element of $R$ at stage of definition $X$. When we work with an arbitrary fixed stage of definition we will sometimes write simply $r \in R$ where it causes no confusion. For interpreting existential quantification and disjunction we will need to consider covers $(\iota_i : X_i \to X)_i$ of the stage of definition $X$. Then if $a \in_X R$ will write $a|_{X_i}$ for the element $a\iota_i \in_{X_i} R$.

Let $D_W$ be a Weil spectrum in $\mathcal{E}$. Then we abuse notation by writing $D_W$ for the object $(M \times D_W, \pi_1)$ of $\mathcal{E}/M$.

**Definition 3.14.** Let $a, b \in_X B$ where $X$ and $B$ are objects of the topos $\mathcal{E}/M$. Then $a \sim b$ iff the proposition

$$\bigvee_{W \in Weil} \exists \phi \in B^{D_W}. \, \exists d \in D_W. \, \phi(0) = a \wedge \phi(d) = b$$

holds in the internal logic of $\mathcal{E}/M$.

Explicitly: there exists a cover $(\iota_i : X_i \to X)_{i \in I}$ in $\mathcal{E}/M$ such that for each $i$ there exists an object $D_{W_i} \in Spec(Weil)$, an arrow $\phi_i : X_i \times D_{W_i} \to B$ and an arrow $d_i : X_i \to D_{W_i}$ such that

$$\begin{array}{ccc} X_i & \xrightarrow{(1_{X_i}, 0)} & X_i \times D_{W_i} \\ {\scriptstyle a|_{X_i}}\downarrow & & \downarrow {\scriptstyle \phi_i} \\ B & \xrightarrow{1_B} & B \end{array}$$

and

$$\begin{array}{ccc} X_i & \xrightarrow{(1_{X_i}, d_i)} & X_i \times D_{W_i} \\ {\scriptstyle b|_{X_i}}\downarrow & & \downarrow {\scriptstyle \phi_i} \\ B & \xrightarrow{1_B} & B \end{array}$$

commute.

**Remark 3.15.** The relation $\sim$ is not always symmetric. In fact it is not symmetric in the case $B = D$ and $M = 1$ as described in Lemma 4.13.

**Definition 3.16.** The relation $\approx$ is the transitive closure of $\sim$ in the internal logic of $\mathcal{E}/M$. This means that for $a, b \in B$ we have $a \approx b$ iff the proposition

$$\bigvee_{n \in \mathbb{N}} \exists \vec{x} \in B^n. \bigwedge_{1 \leq k \leq n-1} (\pi_k \vec{x} \sim \pi_{k+1} \vec{x}) \wedge (\pi_1 \vec{x} = a) \wedge (\pi_n \vec{x} = b)$$

holds in the internal logic of $\mathcal{E}/M$.

In terms of covers: let $a, b \in_X B$ where $X$ and $B$ are objects of $\mathcal{E}/M$. Then $a \approx b$ iff there exists a cover $(\iota_i : X_i \to X)_{i \in I}$ and for each $i$ there exists a natural number $n_i$ and elements $x_{i_0}, x_{i_1}, ..., x_{i_{n_i}} \in_{X_i} B$ such that

$$a|_{X_i} = x_{i_0} \sim x_{i_1} \sim ... \sim x_{i_{n_i}} = b|_{X_i}$$

**Remark 3.17.** For any arrow $f : A \to B$ we have that $a \sim a'$ in $A$ implies that $fa \sim fa'$ in $B$. Indeed if we have $D \in Spec(Weil)$, $\phi \in A^D$ and $d \in D$ such that $\phi(0) = a$ and $\phi(d) = a'$ then for the same $D$ and $d$ we see that $\psi = f^D\phi$ has $\psi(0) = fa$ and $\psi(d) = fa'$.

We can easily iterate this procedure to obtain that $a \approx a'$ in $A$ implies $fa \approx fa'$ in $B$.

**Definition 3.18.** Let $f : A \to B$ be an arrow in $\mathcal{E}/M$. Then $f$ is $\mathcal{W}$-dense (or $f \in L_\mathcal{W}$) iff the proposition

$$\forall b \in B.\ \exists a \in A.\ fa \approx b$$

holds in the internal logic of $\mathcal{E}/M$.

Explicitly: for all $b \in_X B$ there exists a cover $(\iota_i : X_i \to X)_{i \in I}$ and elements $a_i \in_{X_i} A$ such that $f(a_i) \approx b|_{X_i}$.

**Definition 3.19.** Let $g : A \to B$ be an arrow in $\mathcal{E}/M$. Then $g$ is $\mathcal{W}$-closed (or $g \in R_\mathcal{W}$) iff the propositions

$$\forall a \in A.\ \forall b \in B.\ ga \approx b \implies (\exists c \in A.\ gc = b)$$

and

$$\forall a, a' \in A.\ ga = ga' \implies (a \approx a')$$

hold in the internal logic of $\mathcal{E}/M$.

Explicitly the first condition is: for all $a \in_X A$ and $b \in_X B$ such that $ga \approx b$ there exists a cover $(\iota_i : X_i \to X)_{i \in I}$ and elements $c_i \in_{X_i} A$ such that $gc_i = b|_{X_i}$. Since the second condition only uses universal quantification and conjunction it is not necessary to pass to a cover.

**Remark 3.20.** Note that in the sequel the right class of the jet factorisation system will turn out not to be simply $R_\mathcal{W}$ but its intersection with the monomorphisms in $\mathcal{E}/M$. The larger class $R_\mathcal{W}$ will be useful in Section 3.3.

From now on we will work entirely in the internal logic of $\mathcal{E}/M$. The interested reader is welcome to translate the statements below into their external versions involving covers by applying the sheaf semantics explained in Section VI.7 of [14].

**Lemma 3.21.** Let $g : B \rightarrowtail E$ be a $\mathcal{W}$-closed monomorphism. Suppose that $gb \sim gb'$ in $E$. Then $b \sim b'$ in $B$.

*Proof.* Since $gb \sim gb'$ there exists $D \in Spec(Weil)$, $\phi \in E^D$ and $d \in D$ such that $\phi(0) = gb$ and $\phi(d) = gb'$. However it is immediate from the fact that $g$ is $\mathcal{W}$-closed that $\phi$ is in the image of $g^D : B^D \rightarrowtail E^D$ and so there exists $\psi$ such that $\phi = g^D\psi$. But $g(\psi(0)) = gb$ and $g(\psi(d)) = gb'$ hence $\psi(0) = b$ and $\psi(d) = b'$ and $b \sim b'$ as required. □

**Corollary 3.22.** Let $g : B \rightarrowtail E$ be a $\mathcal{W}$-closed monomorphism. Suppose that $gb \approx gb'$ in $E$. Then $b \approx b'$ in $B$.

*Proof.* Let $gb = e_0 \sim e_1 \sim ... \sim e_n = gb'$ exhibit $gb \approx gb'$. Then the fact that $g$ is $\mathcal{W}$-closed combined with $e_0 = gb$ implies that there exists $b_1 \in B$ such that $e_1 = gb_1$. Then by Lemma 3.21 we see that $b \sim b_1$. The result follows easily by iterating this procedure. $\square$

**Lemma 3.23.** *Let $h : A \to E$ be an arrow in $\mathcal{E}/M$. Then there exist $g \in R_\mathcal{W}$ and $f \in L_\mathcal{W}$ such that $g$ is a monomorphism and $h = gf$. The mediating object in the factorisation has the presentation*

$$B = \{x \in E : \exists a \in A. \ ha \approx x\} \xrightarrow{g} E$$

*in the internal logic of $\mathcal{E}/M$.*

*Proof.* It is immediate that $h$ factors through the subobject $B$ because the relation $\approx$ is reflexive. Write $h = gf$ for this factorisation.

To see that $g$ is $\mathcal{W}$-closed let $b \in B$ and $e \in E$ such that $gb \approx e$. By the definition of $B$ there exists an $a \in A$ such that $ha \approx gb$. Hence by the transitivity of $\approx$ we obtain that $ha \approx e$. So $e$ lies in the subobject $B$ and so $g$ is $\mathcal{W}$-closed as required.

To see that $f$ is $\mathcal{W}$-dense let $b \in B$. Now by the definition of $B$ there exists an $a \in A$ such that $ha \approx gb$. But since $g$ is a $\mathcal{W}$-closed monomorphism we can use Corollary 3.22 we deduce that $fa \approx b$ as required. $\square$

**Proposition 3.24.** *Let $\mathcal{M}$ be the class of monomorphisms in $\mathcal{E}/M$. Then the pair*

$$(L, R) = (L_\mathcal{W}, R_\mathcal{W} \cap \mathcal{M})$$

*defines a $(\mathcal{E}/M)$-factorisation system.*

*Proof.* We will check the conditions of Lemma 2.8. The existence of factorisations is Lemma 3.23 and it is clear that the classes $L_\mathcal{W}$ and $R_\mathcal{W} \cap \mathcal{M}$ are replete.

It remains to show that for all $\mathcal{W}$-closed monomorphisms $g : C \rightarrowtail E$ and all $\mathcal{W}$-dense arrows $f : A \to B$ we have that $f \perp_{\mathcal{E}/M} g$. That means we need to show that the square

$$\begin{array}{ccc} C^B & \xrightarrow{C^f} & C^A \\ \downarrow{g^B} & & \downarrow{g^A} \\ E^B & \xrightarrow{E^f} & E^A \end{array}$$

is a pullback. So suppose that $\phi \in E^B$ and $\psi \in C^A$ such that $\phi f = g\psi$. We define $\xi \in C^B$ as follows. Start with $b \in B$. Since $f$ is $\mathcal{W}$-dense there exists $a \in A$ such that $fa \approx b$. Then by Remark 3.17 we have that $g\psi a = \phi fa \approx \phi b$. Now since $g$ is $\mathcal{W}$-closed we have that there exists $c \in C$ such that $gc = \phi b$. This $c$ is unique because $g$ is monic. So finally we define $\xi b = c$. It is immediate that $g\xi b = gc = \phi b$. From the equation $g\xi fa = \phi fa = g\psi a$ we deduce that $\xi fa = \psi a$ as required. $\square$

**Proposition 3.25.** *Let $f : A \rightarrowtail B$ be a monomorphism in $\mathcal{E}/M$. Then $f$ is $\mathcal{W}$-closed iff for all $D \in Spec(Weil)$ the square*

$$
\begin{array}{ccc}
A^{(M \times D, \pi_1)} & \xrightarrow{A^0} & A \\
{\scriptstyle f^D} \downarrow & & \downarrow {\scriptstyle f} \\
B^{(M \times D, \pi_1)} & \xrightarrow{B^0} & B
\end{array}
\tag{3}
$$

*is a pullback.*

*Proof.* We will show that $L_\infty \subset L_\mathcal{W}$ and $R_\infty \subset R_\mathcal{W}$. This will suffice to prove the result because

$$L_\infty \subset L_\mathcal{W} \implies L_\mathcal{W}^\perp \subset L_\infty^\perp \implies R_\mathcal{W} \subset R_\infty$$

To show that $L_\infty \subset L_\mathcal{W}$ we need to show that for all $D \in Spec(Weil)$ the arrow $(M, 1_M) \to (M \times D, \pi_1)$ is in $L_\mathcal{W}$. For this it will suffice to show that for all $b \in (M \times D, \pi_1)$ we have $0 \approx b$. Here $0$ denotes the global element $(1_M, 0) : (M, 1_M) \to (M \times D, \pi_1)$. So we choose $D_W = (M \times D, \pi_1)$, $\phi = 1_{M \times D}$ and $d = b$. Then $\phi(0) = 0$ and $\phi(d) = b$.

To show that $R_\infty \subset R_\mathcal{W}$ let $f$ be a monomorphism, let $a \in A$ and $b \in B$ such that $fa \sim b$ and suppose that the square in (3) is a pullback. The condition $fa \sim b$ means that there is a $D_W \in Spec(Weil)$, a $\phi \in B^{(M \times D_W, \pi_1)}$ and a $d \in (M \times D_W, \pi_1)$ such that $\phi(0) = fa$ and $\phi(d) = b$. Since $\phi(0) = fa$ we can induce a $\psi \in A^{(M \times D, \pi_1)}$ using the pair $(a, \phi)$. But then we have $f\psi(d) = \phi(d) = b$.

We now iterate this argument to obtain that $f$ is $\mathcal{W}$-closed as required. $\square$

**Corollary 3.26.** *The $(L_\mathcal{W}, R_\mathcal{W} \cap \mathcal{M})$ factorisation system and the jet factorisation system coincide in $\mathcal{E}/M$.*

### 3.3 Stability Properties of the Jet Factorisation

Recall that for all factorisation systems the left class is closed under colimits and the right class is closed under limits. The $(Epi, Mono)$-factorisation system has the additional property that the left class is closed under pullbacks. In this section we identify a condition on an arrow $g$ in the left class of the jet factorisation system which guarantees that the pullback of $g$ along a $\mathcal{W}$-closed arrow $k$ is again jet dense.

**Proposition 3.27.** *Let $g$ be jet dense and $k$ be $\mathcal{W}$-closed in $\mathcal{E}/M$. Suppose that the relation $\approx$ is symmetric on the object $E$ and that the square*

$$
\begin{array}{ccc}
A & \xrightarrow{h} & B \\
{\scriptstyle f} \downarrow & & \downarrow {\scriptstyle g} \\
C & \xrightarrow{k} & E
\end{array}
$$

*is a pullback. Then $f$ is also jet dense.*

*Proof.* Recall that an arrow in $\mathcal{E}/M$ is jet dense iff it is $\mathcal{W}$-dense. Let $c \in C$. We need to show that there exists $a \in A$ such that $fa \approx c$. Since $g$ is $\mathcal{W}$-dense there exists $b \in B$ such that $gb \approx kc$. Since $\approx$ is symmetric on $E$ we see that also $kc \approx gb$. Now $k$ is $\mathcal{W}$-closed so there exists $c' \in C$ such that $c \approx c'$ and $kc' = gb$. The $a \in A$ that we require is the one defined by the pair $a = (c', b)$.

We now confirm that $f(c', b) = c' \approx c$. First we see that $kc' = gb \approx kc$ and so there exists $c'' \in C$ such that $c' \approx c''$ and $kc'' = kc$. But now by the definition of $\mathcal{W}$-closed we have that $c'' \approx c$ and by transitivity of $\approx$ that $c' \approx c$ as required. $\square$

**Corollary 3.28.** *Let $g$ be jet dense and $k$ be jet closed. Suppose that the relation $\approx$ is symmetric on $E$ and the square*

$$\begin{array}{ccc} A & \xrightarrowtail{h} & B \\ \downarrow{\scriptstyle f} & & \downarrow{\scriptstyle g} \\ C & \xrightarrowtail{k} & E \end{array}$$

*is a pullback. Then $f$ is also jet dense.*

## 4 The Jet Part Construction

In this section we construct the infinitesimal (or jet) part of an internal category or internal groupoid in a well-adapted model of synthetic differential geometry. For reasons described in Section 4.1 the jet part of a category (or groupoid) will consist of all the arrows of the category that can be reached by a sequence of source constant infinitesimal perturbations from an identity arrow of the category. In Section 4.1 we show that the jet part of a category $\mathbb{C}$ is closed under composition and hence defines a subcategory $\mathbb{C}_\infty$ of $\mathbb{C}$. In Section 4.2 we describe how the subcategory of all internal categories $\mathbb{C}$ for which the inclusion $\iota_\mathbb{C}^\infty$ is an isomorphism is not only a coreflective subcategory of $Cat(\mathcal{E})$ but an $\mathcal{E}$-mono-coreflective subcategory as described in Definition 4.9. Although the jet part $\mathbb{G}_\infty$ of a groupoid $\mathbb{G}$ in $\mathcal{E}$ is therefore a category, it is not necessarily true that $\mathbb{G}_\infty$ is a groupoid. In Section 4.3 we find a necessary and sufficient condition that makes $\mathbb{G}_\infty$ a groupoid: namely that the relation $\approx$ defined in Section 3.2 is symmetric. Then easy to see that the category of jet groupoids satisfying this condition is an $\mathcal{E}$-mono-coreflective subcategory of the category of groupoids satisfying this condition.

### 4.1 The Jet Part of a Category

We define the jet part of a category in a well-adapted mode $\mathcal{E}$ of synthetic differential geometry. Intuitively the arrow space of the jet part will consist of all the elements of the category which we can reach along an infinitesimally small source constant path starting at an identity arrow. We can put the structure of a reflexive graph on these arrows as follows.

**Notation 4.1.** In this section $\mathbb{C}$ will denote a category in $\mathcal{E}$ with underlying reflexive graph

$$\mathbb{C} = \left( C \begin{array}{c} s \\ \leftarrow e \\ \hline t \end{array} M \right)$$

and composition $\mu$. We will write $n_\times C = C \ _t\!\times_s C \ _t\!\times_s \ldots \ _t\!\times_s C$.

**Definition 4.2.** Let

$$M \xrightarrow{e_\infty} C_\infty \xrightarrow{\iota_C^\infty} C$$
$$\begin{array}{c} 1_M \searrow \quad \downarrow s_\infty \quad \swarrow s \\ M \end{array}$$

be the jet factorisation of $e$ in $\mathcal{E}/M$. Then the *jet reflexive graph of* $\mathbb{C}$ is the reflexive graph

$$\mathbb{C}_\infty = \left( C_\infty \begin{array}{c} s_\infty \\ \leftarrow e_\infty \\ \hline t \circ \iota_C^\infty \end{array} M \right)$$

in $\mathcal{E}$.

To equip this reflexive graph $\mathbb{C}_\infty$ with a composition operation we require a slight digression. To understand the reason for this digression we consider the special case that the base space $M = 1$. Then we can make the following straightforward argument. The arrow

$$C_\infty \times M \xrightarrow{1_{C_\infty} \times e_\infty} C_\infty \times C_\infty$$

is jet dense because (as an enriched factorisation system) the left class of the jet factorisation system is closed under products. Then we define the composition on $C_\infty$ to be the unique lift of the following square

$$\begin{array}{ccc} C_\infty \times M & \xrightarrow{\pi_1} & C_\infty \\ 1_{C_\infty} \times e_\infty \downarrow & \nearrow \mu_\infty & \downarrow \iota_C^\infty \\ C_\infty \times C_\infty & \xrightarrow{\mu \circ (\iota_C^\infty \times \iota_C^\infty)} & C \end{array}$$

and the associativity and unit axioms can be seen to hold. However if we now attempt to do the same thing in the slice category $\mathcal{E}/M$ we can still show that the arrow

$$(C_\infty, s_\infty) \xrightarrow{(1_{C_\infty}, e_\infty)} (C_\infty, t_\infty) \times (C_\infty, s_\infty) \cong (C_\infty \ _{t_\infty}\!\times_{s_\infty} C_\infty, t_\infty \pi_1)$$

is jet dense but there is no way to map out of $(C \ _t\!\times_s C, t\pi_1)$ using $\mu$. The problem is that given arrows $f, g \in C$ such that $cod(f) = dom(g)$ the map $t\pi_1$ picks out the 'middle' object $cod(f)$ which cannot be specified from the composite $\mu(f,g)$ alone. We can rescue the idea of using a lift to define the composition by using the results of Section 3.1 to prove that the arrow

$$(C_\infty, s_\infty) \xrightarrow{(1_{C_\infty}, e_\infty t_\infty)} (C_\infty \ _{t_\infty}\!\times_{s_\infty} C_\infty, s_\infty \pi_1)$$

is jet dense in $\mathcal{E}/M$. Then we can proceed in an analogous fashion to the case $M = 1$.

The next result tells us that the map which takes an arrow $g$ of $\mathbb{C}_\infty$ and returns the composable pair $(g, 1_{cod(g)})$ in $2_\times C_\infty$ is jet dense over the source of $g$.

**Lemma 4.3.** *The arrow*

$$C_\infty \xrightarrow{(1_{C_\infty}, e_\infty t_\infty)} 2_\times C_\infty$$
$$\searrow s_\infty \quad \swarrow s_\infty \circ \pi_1$$
$$M$$

*is jet-dense in $\mathcal{E}/M$.*

*Proof.* The arrow

$$M \xrightarrow{e_\infty} C_\infty$$
$$\searrow 1_M \quad \swarrow s_\infty$$
$$M$$

in $\mathcal{E}/M$ is jet dense by the definition of jet part in Definition 4.2. Then by Corollary 3.11 the arrow

$$C_\infty \xrightarrow{(1_{C_\infty}, e_\infty t_\infty)} 2_\times C_\infty$$
$$\searrow 1_{C_\infty} \quad \swarrow \pi_1$$
$$C_\infty$$

obtained by pulling back along $t_\infty$ is jet dense in $\mathcal{E}/C_\infty$. But now by Corollary 3.10 the arrow

$$C_\infty \xrightarrow{(1_{C_\infty}, e_\infty t_\infty)} 2_\times C_\infty$$
$$\searrow s_\infty \quad \swarrow s_\infty \pi_1$$
$$M$$

obtained by postcomposition by $s_\infty$ is jet dense in $\mathcal{E}/M$ as required. $\square$

Now we are in a position to define a composition on the jet part of a category.

**Corollary 4.4.** *Let $\mathbb{C}$ be a category with composition $\mu : C\ {}_t\!\times_s C \to C$. Let $\mathbb{C}_\infty$ be the jet reflexive graph of of $\mathbb{C}$. Then we can make $\mathbb{C}_\infty$ into a category by defining the composition $\mu_\infty : C_\infty\ {}_t\!\times_s C_\infty \to C_\infty$ as the diagonal lift of the following diagram:*

$$\begin{array}{ccc}
(C_\infty, s_\infty) & \xrightarrow{1_{C_\infty}} & (C_\infty, s_\infty) \\
{\scriptstyle (1_{C_\infty}, e_\infty t_\infty)}\downarrow & {\scriptstyle \mu_\infty}\nearrow & \downarrow {\scriptstyle \iota_C^\infty} \\
(2_\times C_\infty, s_\infty \circ \pi_1) & \xrightarrow{\mu \circ (2_\times \iota_C^\infty)} & (C, s)
\end{array}$$

where $\iota_C^\infty$ is jet closed by the definition of $C_\infty$ and $(1_{C_\infty}, e_\infty t_\infty)$ is jet dense by Lemma 4.3. We call the category $\mathbb{C}_\infty$ the *jet part* of $\mathbb{C}$.

*Proof.* The associativity of $\mu_\infty$ is inherited from the associativity of $\mu$. To see this consider the diagram:

$$\begin{array}{ccc}
(3_\times C, s\pi_1) & \xrightarrow{(1_C, \mu)} & (2_\times C, s\pi_1) \\
& \searrow{}^{3_\times \iota_C^\infty} \quad {}^{2_\times \iota_C^\infty}\nearrow & \\
& (3_\times C_\infty, s_\infty \pi_1) \xrightarrow{(1_{C_\infty}, \mu_\infty)} (2_\times C_\infty, s_\infty \pi_1) & \\
\downarrow{}^{(\mu, 1_C)} & \downarrow{}^{(\mu_\infty, 1_{C_\infty})} \quad \downarrow{}^{\mu_\infty} & \downarrow{}^{\mu} \\
& (2_\times C_\infty, s_\infty \pi_1) \xrightarrow{\mu_\infty} (C_\infty, s_\infty) & \\
& \nearrow{}^{2_\times \iota_C^\infty} \quad {}^{\iota_C^\infty}\searrow & \\
(2_\times C, s\pi_1) & \xrightarrow{\mu} & (C, s)
\end{array}$$

where the outer square commutes because $\mu$ is associative and the top, bottom, left and right squares commute using the definition of $\mu_\infty$ above. But this implies that the inner square commutes because $\iota_C^\infty$ is monic.

One of the unit laws for $\mu_\infty$ is already enforced by the upper commutative triangle in the definition of $\mu_\infty$. The other follows from combining the fact that $\iota_C^\infty$ is monic and that in the diagram:

$$\begin{array}{ccc}
(C, s) & \xrightarrow{1_C} & (C, s) \\
& \searrow{}^{\iota_C^\infty} \quad {}^{\iota_C^\infty}\nearrow & \\
& (C_\infty, s_\infty) \xrightarrow{1_{C_\infty}} (C_\infty, s_\infty) & \\
\downarrow{}^{(es, 1_C)} & \downarrow{}^{(e_\infty s_\infty, 1_{C_\infty})} \quad \downarrow{}^{1_{C_\infty}} & \downarrow{}^{1_C} \\
& (2_\times C_\infty, s_\infty \pi_1) \xrightarrow{\mu_\infty} (C_\infty, s_\infty) & \\
& \nearrow{}^{2_\times \iota_C^\infty} \quad {}^{\iota_C^\infty}\searrow & \\
(2_\times C, s\pi_1) & \xrightarrow{\mu} & (C, s)
\end{array}$$

the outer square commutes using a unit law for $\mu$ and the other squares are immediately seen to commute. $\square$

### 4.2 The Category of Jet Categories

**Definition 4.5.** Let $\mathbb{C}$ be a category in $\mathcal{E}$ and $\mathbb{C}_\infty$ be the category on its jet part as defined in Corollary 4.4. Then $\mathbb{C}$ is a *jet category* iff the inclusion $\iota_\mathbb{C}^\infty : \mathbb{C}_\infty \rightarrowtail \mathbb{C}$ induced by $\iota_\infty$ is an isomorphism. We write $Cat_\infty(\mathcal{E})$ for the full subcategory of $Cat(\mathcal{E})$ on the jet categories.

**Lemma 4.6.** *The function* $(-)_\infty : Cat(\mathcal{E}) \to Cat_\infty(\mathcal{E})$ *extends to a functor.*

*Proof.* Let $\phi : \mathbb{C} \to \mathbb{B}$ be an internal functor. Then the square

$$\begin{array}{ccc} (M, 1_M) & \xrightarrow{e_\infty^B \phi_0} & (B_\infty, s_\infty^B) \\ {\scriptstyle e_\infty^C} \downarrow & \phi_\infty \nearrow & \downarrow {\scriptstyle \iota_\infty^B} \\ (C_\infty, s_\infty^C) & \xrightarrow{\phi_1 \iota_\infty^C} & (B, s^B) \end{array}$$

commutes in $\mathcal{E}$ and hence there exists a unique filler $\phi_\infty$. It is immediate from the definition that $\phi_\infty$ preserves identities. We now remark that in the cube

$$\begin{array}{ccccc} (2_\times C, s^C \pi_1) & & \xrightarrow{\mu^C} & & (C, s^C) \\ & \nwarrow {\scriptstyle 2_\times \iota_\infty^C} & & {\scriptstyle \iota_\infty^C} \nearrow & \\ & & (2_\times C_\infty, s_\infty^C \pi_1) \xrightarrow{\mu_\infty^C} (C_\infty, s_\infty^C) & & \\ {\scriptstyle 2_\times \phi_1} \downarrow & & {\scriptstyle 2_\times \phi_\infty} \downarrow \quad \downarrow {\scriptstyle \phi_\infty} & & \downarrow {\scriptstyle \phi_1} \\ & & (2_\times B_\infty, s_\infty^B \pi_1) \xrightarrow{\mu_\infty^B} (B_\infty, s^B) & & \\ & \swarrow {\scriptstyle 2_\times \iota_\infty^B} & & {\scriptstyle \iota_\infty^B} \searrow & \\ (2_\times B, s^B \pi_1) & & \xrightarrow{\mu^B} & & (B, s^B) \end{array}$$

the outer square commutes by functoriality of $\phi$, the left and right faces commute by definition of $\phi_\infty$ and the top and bottom faces commute by the definition of $\mu_\infty$. Therefore the inner square commutes because the arrow $\iota_\infty^B$ is a monomorphism and hence $\phi_\infty$ preserves composition as required. $\square$

**Proposition 4.7.** *We have an adjunction $j \dashv (-)_\infty$ where $j$ is the full inclusion $Cat_\infty(\mathcal{E}) \hookrightarrow Cat(\mathcal{E})$. In other words $Cat_\infty(\mathcal{E})$ is a coreflective subcategory of $Cat(\mathcal{E})$.*

*Proof.* Let $\mathbb{K}$ be a jet category; this means that the inclusion $\iota_\infty^K : \mathbb{K}_\infty \rightarrowtail \mathbb{K}$ is an isomorphism. We define the unit $\eta$ by $\eta_\mathbb{K} = (\iota_\infty^K)^{-1}$. Let $\mathbb{C}$ be an arbitrary category in $\mathcal{E}$. We define the counit $\varepsilon$ of the adjunction by $\varepsilon_\mathbb{C} = \iota_\infty^C$. Then $\varepsilon_{j(\mathbb{K})} \circ j(\eta_\mathbb{K}) = \iota_\infty^K \circ (\iota_\infty^K)^{-1} = 1_{j\mathbb{K}}$ and $(\varepsilon_\mathbb{C})_\infty \circ \eta_{\mathbb{C}_\infty} = (\iota_\infty^C)_\infty \circ (\iota_\infty^{C_\infty})^{-1}$. But by definition of $(\iota_\infty^C)_\infty$ we see that

$$\begin{array}{ccccc} M & \xrightarrowtail{e_\infty^{C_\infty}} & (\mathbb{C}_\infty)_\infty & \xrightarrowtail{\iota_\infty^{C_\infty}} & \mathbb{C}_\infty \\ \downarrow {\scriptstyle 1_M} & & \downarrow {\scriptstyle (\iota_\infty^C)_\infty} & & \downarrow {\scriptstyle \iota_\infty^C} \\ M & \xrightarrowtail{e_\infty^C} & \mathbb{C}_\infty & \xrightarrowtail{\iota_\infty^C} & \mathbb{C} \end{array}$$

commutes and so $\iota_\infty^C \circ (\iota_\infty^C)_\infty \circ (\iota_\infty^{C_\infty})^{-1} = \iota_\infty^C$. Hence $(\iota_\infty^C)_\infty \circ (\iota_\infty^{C_\infty})^{-1} = 1_{\mathbb{C}_\infty}$ because $\iota_\infty^C$ is a monomorphism. $\square$

In fact $Cat_\infty(\mathcal{E})$ is an $\mathcal{E}$-mono-coreflective subcategory as described in the introduction.

**Lemma 4.8.** *Let $\mathbb{C}$ and $\mathbb{B}$ be internal categories in $\mathcal{E}$ and $\mathbb{B}_\infty$ and $\mathbb{C}_\infty$ be their jet parts. Then $\mathbb{C}^{\mathbb{B}_\infty} \cong \mathbb{C}_\infty^{\mathbb{B}_\infty}$ in $\mathcal{E}$.*

*Proof.* To show that $\mathbb{C}^{\mathbb{B}_\infty} \cong \mathbb{C}_\infty^{\mathbb{B}_\infty}$ it will suffice to show that for all representable objects $X$ in $\mathcal{E}$ and internal functors $F : \mathbb{B}_\infty \times \dot{X} \to \mathbb{C}$ we have a unique lift $G$ making

$$\begin{array}{ccc} & & \mathbb{C}_\infty \\ & \nearrow{G} & \downarrow{\iota_\infty^C} \\ \mathbb{B}_\infty \times \dot{X} & \xrightarrow{F} & \mathbb{C} \end{array}$$

commute. But we can just take $G = F_\infty$ because the fact that $(-)_\infty$ is a right adjoint implies that $(\mathbb{B}_\infty \times \dot{X})_\infty = \mathbb{B}_\infty \times \dot{X}$. $\square$

This means that $Cat_\infty(\mathcal{E})$ is an $\mathcal{E}$-mono-coreflective subcategory of $Cat(\mathcal{E})$ defined as follows:

**Definition 4.9.** Let $\mathcal{U}$ and $\mathcal{C}$ be $\mathcal{E}$-categories with underlying categories $\mathcal{U}_0$ and $\mathcal{C}_0$ respectively. Then the category $\mathcal{U}$ is an *$\mathcal{E}$-mono-coreflective subcategory* of $\mathcal{C}$ iff there is an adjunction

$$\mathcal{U}_0 \xrightarrow[\underset{(-)_\infty}{\longleftarrow}]{\longrightarrow} \mathcal{C}_0$$

such that the left adjoint is full and faithful, the counit $\iota^\infty$ is a monomorphism and for all objects $\mathbb{U}$ of $\mathcal{U}$ and $\mathbb{C}$ of $\mathcal{C}$ the arrow

$$\mathbb{C}_\infty^{\mathbb{U}} \xrightarrow{(\iota_\mathbb{C}^\infty)^{\mathbb{U}}} \mathbb{C}^{\mathbb{U}}$$

is an isomorphism in $\mathcal{E}$.

### 4.3 The Jet Part of a Groupoid

In order to put the structure of a groupoid on the jet part $\mathbb{C}_\infty$ of a category $\mathbb{C}$ we require a little more than the assumption that $\mathbb{C}$ has the necessary additional structure and relations to make it a groupoid. For the rest of this section we fix a groupoid $\mathbb{G}$ that has underlying reflexive graph

$$G \xrightarrow[\underset{t}{\longrightarrow}]{\overset{s}{\longrightarrow}} M$$
$$\xleftarrow{e}$$

and multiplication $\mu$. We identify a necessary condition for the jet part $\mathbb{G}_\infty$ of $\mathbb{G}$ to have groupoid structure. First we need a preparatory lemma.

**Lemma 4.10.** *Let $a, b \in_{(X,sa)} (G, s)$ be generalised elements in $\mathcal{E}/M$ such that $a \approx b$ at stage of definition $(X, sa)$. Let $c \in_{(X,sc)} (G, s)$ such that $tc = sa(= sb)$. Then $ac \approx bc$ at stage of definition $(X, sc)$.*

*Proof.* Suppose that $a \approx b$ is witnessed by the following data:

- a cover $(\iota_i : (X_i, sa_i) \to (X, sa))_{i \in I}$ and for each $i \in I$:
    - an arrow $\phi_i : (X_i \times D_i, sa_i \pi_1) \to (G, s)$;
    - an arrow $d_i : (X_i, sa_i) \to (M \times D_i, \pi_1)$;

such that $\phi(1_{X_i}, 0) = a_i$ and $\phi_i(1_{X_i}, \pi_2 d_i) = b_i$ where $a_i$ and $b_i$ are the restrictions of $a$ and $b$ respectively to $X_i$.

As a first step we show that $(c, a) \approx (c, b)$ as generalised elements at stage of definition $(X_i, sc_i)$ where $c_i$ is the restriction of $c$ to $X_i$. To do this we choose:

- the cover $(\iota_i : (X_i, sc_i) \to (X, sc))_{i \in I}$ and for each $i \in I$:
    - the arrow $\overline{\phi_i} = (c_i \pi_1, \phi_i) : (X_i \times D_i, sc_i \pi_1) \to (2_\times G, s\pi_1)$;
    - the arrow $\overline{d_i} = (sc_i, \pi_2 d_i) : (X_i, sc_i) \to (m \times D_i, \pi_1)$;

and note that:
$$(c_i \pi_1, \phi_i)(1_{X_i}, 0) = (c_i, a_i)$$
and
$$(c_i \pi_1, \phi_i)(1_{X_i}, \pi_2 d_i) = (c_i, b_i)$$
hold. Hence $(c, a) \approx (c, b)$. But now the result follows from Lemma 3.17 by applying $\mu$. $\square$

**Proposition 4.11.** *Let $\mathbb{G}$ be a groupoid in $\mathcal{E}$ with arrow space $G$ and object space $M$. Suppose further that $\mathbb{G}_\infty$ is a groupoid. Then the relation $\approx$ is symmetric on $(G, s)$ in $\mathcal{E}/M$.*

*Proof.* Let $a, b \in_{(X, sa)} (G, s)$ such that $a \approx b$ at stage of definition $(X, sa)$. Then $a^{-1} \in_{(X, ta)} (G, s)$ has $ta^{-1} = sa (= sb)$. So by precomposing with $a^{-1}$ and using Lemma 4.10 we have that $eta = aa^{-1} \approx ba^{-1}$ at stage of definition $(X, ta)$ and hence $ba^{-1} \in_{(X, ta)} (G_\infty, s_\infty)$. Since $\mathbb{G}_\infty$ is a groupoid we have that $ab^{-1} \in_{(X, ta)} (G_\infty, t_\infty)$ also and hence $ab^{-1} \in_{(X, tb)} (G_\infty, s_\infty)$.

This means that $etb \approx ab^{-1}$ at stage of definition $(X, tb)$. Now we note that $b \in_{(X, sb)} (G, s)$ has $tb = setb = sab^{-1}$ and so by Lemma 4.10 again we deduce that $b \approx ab^{-1}b = a$ as required. $\square$

Now we give a counterexample which shows that we cannot immediately specialise the jet part construction for categories to construct a jet part for an arbitrary groupoid in $\mathcal{E}$. We will use one of the simplest non-classical groupoids we have at our disposal: the pair groupoid $\nabla D$ where $D = \{x \in R : x^2 = 0\}$.

**Lemma 4.12.** *The jet part of the pair groupoid on the object $D$ has the following arrow space:*
$$(\nabla D)^2_\infty = \{(a, b) \in D \times D : a \approx_1 b\}$$
*where $\approx_1$ denotes the neighbour relation in $\mathcal{E}/1$.*

*Proof.* Recall from Lemma 3.23 that the arrow space of $\mathbb{I}_\infty$ is characterised as follows. An arrow $(a,b) : (X, \xi) \to (I^2, \pi_1)$ in $\mathcal{E}/I$ factors through $(\mathbb{I}^2_\infty, \pi_1)$ iff there exists a cover $(\iota_i : (X_i, \xi_i) \to (X, \xi))_i$ such that for all $i$ there exist $W_i \in Spec(Weil)$, $\phi_i : (X_i, \xi_i) \times (I \times D_{W_i}, \pi_1) \to (I^2, \pi_1)$ and $d_i : (X_i, \xi_i) \to (I \times D_{W_i}, \pi_1)$ making

$$\begin{array}{c}
(X_i, \xi_i) \\
{\scriptstyle (1_{X_i}, d_i)} \downarrow \quad \searrow {\scriptstyle (a_i, b_i)} \\
(X_i \times D_{W_i}, \xi_i \pi_1) \xrightarrow{\phi_i} (I^2, \pi_1) \\
{\scriptstyle (1_{X_i}, 0)} \uparrow \quad \nearrow {\scriptstyle (m_i, m_i)} \\
(X_i, \xi_i)
\end{array}$$

commute where $a_i$, $b_i$, $\phi_i$ and $m_i$ are the restrictions of $a$, $b$, $\phi$ and $m$ to $X_i$. Hence $m_i = a_i$ and $(a,b)$ factors through $\mathbb{I}^2_\infty$ iff $a \approx b$ in $\mathcal{E}/1$. $\square$

Therefore to show that $(\nabla D)^2_\infty$ is not a groupoid it will suffice to show that $\approx$ is not symmetric on $D$ in $\mathcal{E}/1$. To prove this we will show that any jet starting from the generalised element $1_D$ must be trivial. The intuitive reason for this is that $D$ is not closed under addition and so there is no more 'space' for the jet to move into.

**Lemma 4.13.** *The relation $\approx$ is not symmetric on $D$.*

*Proof.* Let us consider the generalised elements at stage $D$ described by $0 : D \to D$ and $1_D$. It will suffice to show that $0 \sim 1_D$ but not $1_D \approx 0$. To see that $0 \sim 1_D$ we choose $D_W = D$, $\phi = 1_D$ and $d = 1_D$. Then $\phi(0) = 0$ and $\phi(d) = 1_D$.

To show that $1_D \approx 0$ does not hold it will suffice to show that for all elements $f$ such that $1_D \sim f$ then necessarily $f = 1_D$. So let us suppose that we have an $f$ such that $1_D \sim f$. Since the only covers of $D$ are trivial this would mean that there exist $D_W \in Spec(Weil)$, $\phi : D \times D_W \to D$ and $d : D \to D_W$ such that $\phi(x, 0) = x$ and $\phi(x, d(x)) = f(x)$ for all $x \in D$. Let $w$ be the number of indeterminates in the polynomial defining the Weil presentation $W$. Now we use Hadamard's Lemma twice and the fact that $D$ is defined by the formula $x^2 = 0$ to see that

$$\phi(x_1, \vec{x}) \cong \phi_0(\vec{x}) + x_1 \phi_1(\vec{x})$$

for some smooth functions $\phi_0, \phi_1 : \mathbb{R}^w \to \mathbb{R}$. Now the equation $\phi(a, 0) = a$ tells us that

$$\phi_0(0) + x_1 \phi_1(0) = \phi(x_1, 0) = x_1$$

and so $\phi_0(0) = 0$ and $\phi_1(0) = 1$. Hence by Hadamard's Lemma we see that

$$\phi_1(\vec{x}) = 1 + \Sigma_{i=2}^{w+1} x_i \psi_i(\vec{x})$$

for some $\psi_i : \mathbb{R}^w \to \mathbb{R}$. But since for all $i$ there is an equality of the form $x_i^{k_i} = 0$ in $W$ we see that $N = \Sigma_{i=2}^{w+1} x_i \psi_i(\vec{x})$ is nilpotent of degree $n = \Sigma_{i=2}^{w+1} k_i$. (This follows from the pigeonhole principle.) Therefore the arrow

$$i_\phi = \Sigma_{j=0}^{n-1}(-1)^j N^k : D_W \to D$$

is a pointwise multiplicative inverse for $\phi_1$. Now because $\phi$ has codomain $D$ we must have that

$$\phi_0(\vec{x})^2 + 2x_1\phi_0(\vec{x})\phi_1(\vec{x}) = \phi(x_1,\vec{x})^2 = 0$$

and so $\phi_0(\vec{x})\phi_1(\vec{x}) = 0$. But since $\phi_1$ has a pointwise multiplicative inverse this means that $\phi_0(\vec{x}) = 0$ and so $\phi(x_1,\vec{x}) \cong x_1\phi_1(\vec{x})$. Similarly we see that

$$d(x) \cong \vec{a} + \vec{b}x$$

where $(a_i + b_i x)^{k_i} = 0$ when $x^2 = 0$. But since $a_i \in \mathbb{R}$ we see that $a_i = 0$ and hence

$$\phi(x, d(x)) = x + x\Sigma_{i=2}^{w+1} d(x)_i \psi_i(d(x)) = x + x\Sigma_{i=2}^{w+1} b_i x \psi_i(d(x)) = x$$

and we deduce that $f = 1_D$ as required. $\square$

**Corollary 4.14.** *The jet part $(\nabla D)_\infty$ of the pair groupoid $\nabla D$ is not a groupoid.*

*Proof.* The result follows immediately from Lemma 4.13 and the remarks preceding it. $\square$

Fortunately the condition that the relation $\approx$ is symmetric on $(G, s)$ in $\mathcal{E}/M$ is not only necessary but also sufficient to ensure that the jet part $\mathbb{G}_\infty$ of $\mathbb{G}$ is a groupoid.

**Lemma 4.15.** *Let $a \in (G, s)$ such that $esa \approx a$ in $(G, s)$. Suppose further that $\approx$ is symmetric on $(G, s)$. Then $eta \approx a^{-1}$ in $(G, s)$.*

*Proof.* Since $\approx$ is symmetric we have that $a \approx esa$ and $ta^{-1} = sa$. So by Lemma 4.10 we have that $eta \approx a^{-1}$. $\square$

**Lemma 4.16.** *Let $a, b \in (G, s)$ such that $a \approx b$ in $(G, s)$. Then $a^{-1} \approx b^{-1}$ in $(G, t)$.*

*Proof.* Immediate from Lemma 3.17. $\square$

**Corollary 4.17.** *If $\approx$ is symmetric on $(G, s)$ then the arrow*

$$e_\infty : (M, 1_M) \to (G_\infty, t_\infty)$$

*is jet dense.*

*Proof.* Let $a \in (G_\infty, t_\infty)$. By definition of $G_\infty$ this means that that $esa \approx a$ in $(G, s)$. Since $\approx$ is symmetric on $(G, s)$ we have that $a \approx esa$. Precomposing with $a^{-1}$ and using Lemma 4.15 gives $eta \approx a^{-1}$. Finally applying $(-)^{-1}$ and using Lemma 4.16 gives $eta \approx a$ as required. $\square$

**Proposition 4.18.** *Let $\mathbb{G}$ be a groupoid in $\mathcal{E}$ such that the relation $\approx$ is symmetric on the object $(G_\infty, s_\infty)$ in $\mathcal{E}/M$. Then the jet part $\mathbb{G}_\infty$ can be given the structure of a groupoid.*

*Proof.* By Corollary 4.17 we see that the left arrow in the square

$$\begin{array}{ccc} (M, 1_M) & \xrightarrow{e_\infty} & (G_\infty, s_\infty) \\ \downarrow e_\infty \quad i_{G_\infty} \nearrow & & \downarrow \iota_G^\infty \\ (G_\infty, t_\infty) & \xrightarrow{i_G \iota_G^\infty} & (G, s) \end{array}$$

is jet dense. This means that there is an unique filler $i_{G_\infty}$ which we will take as the inverse for the jet part $\mathbb{G}_\infty$. Since the equations $s_\infty i_{G_\infty} = t_\infty$ and $t_\infty i_{G_\infty} = s_\infty$ are immediately seen to hold it remains to check that the inverse axioms hold. So observe that in the diagram:

$$\begin{array}{ccccc} (G_\infty, s_\infty) & \xrightarrow{(1_{G_\infty}, i_{G_\infty})} & (2_* G_\infty, s_\infty \pi_1) & \xrightarrow{\mu_\infty} & (G_\infty, s_\infty) \\ \downarrow \iota_G^\infty & & \downarrow 2_* \iota_G^\infty & & \downarrow \iota_G^\infty \\ (G, s) & \xrightarrow{(1_G, i_G)} & (2_* G, s\pi_1) & \xrightarrow{\mu} & (G, s) \end{array}$$

the right-hand square commutes by the definition of $\mu_\infty$ in Definition 4.4 and the left-hand square commutes by the definition of $i_{G_\infty}$ above. But now we notice that the bottom row is equal to $1_G$ because $i_G$ is an inverse for the multiplication $\mu$; hence the top row is equal to $1_{G_\infty}$ because $\iota_G^\infty$ is monic. Similarly the diagram

$$\begin{array}{ccccc} (G_\infty, s_\infty) & \xrightarrow{(i_{G_\infty}, 1_{G_\infty})} & (2_* G_\infty, s_\infty \pi_1) & \xrightarrow{\mu_\infty} & (G_\infty, s_\infty) \\ \downarrow \iota_G^\infty & & \downarrow 2_* \iota_G^\infty & & \downarrow \iota_G^\infty \\ (G, s) & \xrightarrow{(i_G, 1_G)} & (2_* G, s\pi_1) & \xrightarrow{\mu} & (G, s) \end{array}$$

shows that the other inverse axiom holds. $\square$

Now we record the analogous result to Lemma 4.8 and Proposition 4.7 when working with groupoids rather than categories.

**Corollary 4.19.** *Let $Gpd^{sym}(\mathcal{E})$ be the full subcategory of $Gpd(\mathcal{E})$ consisting of the groupoids $\mathbb{G}$ for which $\approx$ is symmetric on the arrow space $\mathbb{G}^2$. Let $Gpd_\infty^{sym}(\mathcal{E})$ be the full subcategory of $Gpd^{sym}(\mathcal{E})$ consisting of those groupoids $\mathbb{G}$ for which the subgroupoid inclusion $\iota_\mathbb{G}^\infty : \mathbb{G}_\infty \rightarrowtail \mathbb{G}$ is an isomorphism. Then $Gpd_\infty^{sym}(\mathcal{E})$ is an $\mathcal{E}$-mono-coreflective subcategory of $Gpd^{sym}(\mathcal{E})$.*

## 5  Axiomatics for Lie's Second Theorem

In the remainder of the paper we assume the existence of three pieces of data and prove Lie's second theorem relative to them. The three pieces of data are:

- a topos $\mathcal{E}$,
- an $\mathcal{E}$-mono-coreflective subcategory $Cat_\infty(\mathcal{E})$ of $Cat(\mathcal{E})$,
- an arrow $l : \mathbf{2} \to \mathbb{I}$ in $Cat(\mathcal{E})$.

The symbol $\mathbf{2}$ denotes the category on two objects with a single non-trivial arrow and $\mathbb{I}$ is an arbitrary category in $Cat(\mathcal{E})$. The main intended application is when $\mathcal{E}$ is a well-adapted model of synthetic differential geometry. In this case the adjunction

$$Cat_\infty(\mathcal{E}) \underset{(-)_\infty}{\overset{i}{\rightleftarrows}} Cat(\mathcal{E})$$

is an $\mathcal{E}$-mono-coreflective subcategory by Lemma 4.8 and Proposition 4.7. For the object $\mathbb{I}$ we use the category that has underlying reflexive graph

$$\{(a,b) \in I^2 : a \leq b\} \underset{\pi_2}{\overset{\pi_1}{\underset{\leftarrow \Delta \rightarrow}{\rightrightarrows}}} I$$

and the only possible composition. We have written $I$ for the unit interval in $\mathcal{E}$. The arrow $l$ is then given by $(0,1) \in I^2$.

Note that by working at this greater generality our results will easily adapt to the case in which we replace categories in $\mathcal{E}$ with groupoids in $\mathcal{E}$. This is because $Gpd_\infty(\mathcal{E})$ is a $\mathcal{E}$-mono-coreflective subcategory of $Gpd(\mathcal{E})$ by Corollary 4.19 and we can choose $\mathbb{I}$ to be the pair groupoid on the unit interval $I$. The arrow $l$ is then given by $l = (0,1) \in I^2$ as before.

In Section 6 we use the data $\mathcal{E}$, $(-)_\infty$ and $l$ to define a $Cat(\mathcal{E})$-factorisation system on $Cat(\mathcal{E})$ called the integral factorisation system and generate from it a reflective subcategory

$$Cat(\mathcal{E}) \underset{j}{\overset{(-)_{int}}{\rightleftarrows}} Cat_{int}(\mathcal{E})$$

In Section 7 we work out the appropriate enriched analogues of the connectedness conditions required in classical Lie theory. Then in Section 8 we prove that when we assert these connectedness conditions the functor $(-)_\infty$ is full and faithful.

## 6 Integration of Infinitesimal Paths

As a first step to proving Lie's second theorem we will describe a class of internal categories in $\mathcal{E}$ for which we can integrate paths of infinitesimals ($A$-paths) to macroscopic paths ($G$-paths). We use the theory of enriched factorisation systems to pick out this class.

## 6.1 The Integral Factorisation System

In this section we create an enriched factorisation system which captures the idea of integrating $A$-paths to $G$-paths.

**Definition 6.1.** Let $\mathcal{E}$, $(-)_\infty$ and $\mathbb{I}$ be as in Section 5 and set $\mathcal{E} = \mathcal{V} = Cat(\mathcal{E})$. The *integral factorisation system* is the $Cat(\mathcal{E})$-factorisation system generated by the singleton set

$$\Sigma = \{\mathbb{I}_\infty \xrightarrow{\iota_{\mathbb{I}}^\infty} \mathbb{I}\}$$

using Proposition 2.10.

**Remark 6.2.** Explicitly, an arrow $r$ is in the right class of the integral factorisation system iff

$$\begin{array}{ccc} \mathbb{X}^{\mathbb{I}} & \xrightarrow{\mathbb{X}^{\iota_{\mathbb{I}}^\infty}} & \mathbb{X}^{\mathbb{I}_\infty} \\ \downarrow r^{\mathbb{I}} & & \downarrow r^{\mathbb{I}_\infty} \\ \mathbb{Y}^{\mathbb{I}} & \xrightarrow{\mathbb{Y}^{\iota_{\mathbb{I}}^\infty}} & \mathbb{Y}^{\mathbb{I}_\infty} \end{array}$$

is a pullback in $Cat(\mathcal{E})$ and an arrow $l$ is in the left class of the integral factorisation system iff for all $r$ in the right class

$$\begin{array}{ccc} \mathbb{X}^{\mathbb{B}} & \xrightarrow{\mathbb{X}^l} & \mathbb{X}^{\mathbb{B}} \\ \downarrow r^{\mathbb{B}} & & \downarrow r^{\mathbb{A}} \\ \mathbb{Y}^{\mathbb{B}} & \xrightarrow{\mathbb{Y}^l} & \mathbb{Y}^{\mathbb{A}} \end{array}$$

is a pullback in $Cat(\mathcal{E})$.

**Remark 6.3.** By construction the arrow $(\iota_{\mathbb{I}}^\infty)^n : \mathbb{I}_\infty^n \to \mathbb{I}^n$ is in the left class of the integral factorisation system for all $n \in \mathbb{N}$ and so

$$\begin{array}{ccc} \mathbb{X}^{\mathbb{I}^n} & \xrightarrow{\mathbb{X}^{(\iota_{\mathbb{I}}^\infty)^n}} & \mathbb{X}^{\mathbb{I}_\infty^n} \\ \downarrow r^{\mathbb{I}^n} & & \downarrow r^{\mathbb{I}_\infty^n} \\ \mathbb{Y}^{\mathbb{I}^n} & \xrightarrow{\mathbb{Y}^{(\iota_{\mathbb{I}}^\infty)^n}} & \mathbb{Y}^{\mathbb{I}_\infty^n} \end{array}$$

is a pullback for all $r$ in the right class of the integral factorisation system. Note that the arrow $\iota_{\partial \mathbb{I}^2}^\infty : \partial \mathbb{I}_\infty^2 \to \partial \mathbb{I}^2$ is not in general in the left class of the integral factorisation system. This justifies the use of the simply connectedness condition in Lemma 7.3 and using two dimensional data in general.

## 6.2 The Category of Integral Complete Categories

In this section we recall that the integral factorisation system generates a reflective subcategory of $Cat(\mathcal{E})$.

**Definition 6.4.** The category $Cat_{int}(\mathcal{E})$ is the full subcategory of $Cat(\mathcal{E})$ whose objects are those categories $\mathbb{C}$ for which the arrow

$$\mathbb{C}^{\mathbb{I}} \xrightarrow{\iota_{\mathbb{C}}^{\infty}} \mathbb{C}^{\mathbb{I}_{\infty}}$$

is an isomorphism in $Cat(\mathcal{E})$.

Using the relationship between factorisation systems, completion operations and reflective subcategories in [10] we see that the category $Cat_{int}(\mathcal{E})$ is a reflective subcategory

$$Cat(\mathcal{E}) \underset{j}{\overset{(-)_{int}}{\rightleftarrows}} Cat_{int}(\mathcal{E})$$

of $Cat(\mathcal{E})$ with reflector $(-)_{int}$ that takes a category $\mathbb{C}$ to the mediating object of the integral factorisation of the arrow $! : \mathbb{C} \to 1$. Combining this adjunction with the coreflection $(-)_{\infty}$ gives an adjunction

$$Cat_{\infty}(\mathcal{E}) \underset{(-)_{\infty}}{\overset{(-)_{int}}{\rightleftarrows}} Cat_{int}(\mathcal{E})$$

which is analogous to the adjunction between the category of Lie groups and the category of formal group laws.

# 7 Connectedness, Path Spaces and Global Properties

The lifting property at the core of Lie's second theorem involves lifting internal functors $\mathbb{C}_{\infty} \to \mathbb{X}$ to functors $\mathbb{C} \to \mathbb{X}$. The second stage in our proof of Lie's second theorem is to reformulate this lifting property in terms of generalised elements at stage of definition $\mathbb{I}^n$ where $n \in \{0, 1, 2\}$. It turns out that we need to keep track of the boundary and degeneracy maps between these generalised elements of $\mathbb{C}$ and so in Section 7.2 we organise them into a 2-truncated cubical object in the topos $\mathcal{E}$. In Section 7.1 we define the notions of enriched path and simply connected category and in Section 7.3 and Section 7.4 we show how to describe functors out of a simply connected internal category $\mathbb{C}$ in $\mathcal{E}$ in terms of truncated cubical objects. Then in Section 7.5 we reformulate the lifting problem at the core of Lie's second theorem in terms of truncated cubical objects.

## 7.1 Enriched Connectedness

In classical Lie theory we study how much of the data in a Lie groupoid can be recovered from the subset of this data that is infinitely close to the identity arrows of the Lie groupoid. Since global features such as connectedness cannot be captured by infinitesimal arrows we need to restrict our attention to Lie groupoids that are source simply connected.

We say that a Lie groupoid $\mathbb{G}$ is *source path/source simply connected* iff all of its source fibres are path/simply connected. Let $\nabla I$ be the pair groupoid on the unit interval $I$ that has underlying reflexive graph:

$$I \times I \underset{\pi_2}{\overset{\pi_1}{\underset{\leftarrow \Delta}{\rightrightarrows}}} I$$

with the only possible composition. Then it is easy to see that groupoid homomorphisms $\nabla I \to \mathbb{G}$ are equivalent to arrows $I \to G$ that are source constant and start at an identity element of $G$. Therefore $\mathbb{G}$ is source path connected iff

$$\Gamma(\mathbb{G}^{\nabla I}) \xrightarrow{\Gamma(\mathbb{G}^{\iota \nabla I})} \Gamma(\mathbb{G}^{\partial \nabla I})$$

is an epimorphism in *Set*. We have written $\Gamma$ for the global sections functor and $\partial \nabla I$ for the pair groupoid on the boundary of $I$. Similarly $\mathbb{G}$ is source simply connected iff it is source path connected and

$$\Gamma(\mathbb{G}^{\nabla I^2}) \xrightarrow{\Gamma(\mathbb{G}^{\iota \nabla I^2})} \Gamma(\mathbb{G}^{\partial \nabla I^2})$$

is an epimorphism in *Set*. We have written $\partial \nabla I^2$ for the pair groupoid on the boundary of $I^2$.

When we work with arbitrary groupoids in $\mathcal{E}$ it is more natural to work with epimorphisms between objects of $\mathcal{E}$ than between their sets of global sections. Hence we make the following definitions:

**Definition 7.1.** A groupoid $\mathbb{G}$ in $\mathcal{E}$ is $\mathcal{E}$-*path connected* iff

$$\mathbb{G}^{\nabla I} \xrightarrow{\mathbb{G}^{\iota \nabla I}} \mathbb{G}^{\partial \nabla I}$$

is an epimorphism in $\mathcal{E}$. A groupoid $\mathbb{G}$ in $\mathcal{E}$ is $\mathcal{E}$-*simply connected* iff it is $\mathcal{E}$-path connected and

$$\mathbb{G}^{\nabla I^2} \xrightarrow{\mathbb{G}^{\iota \nabla I^2}} \mathbb{G}^{\partial \nabla I^2}$$

is an epimorphism in $\mathcal{E}$.

This means that for an arbitrary groupoid in $\mathcal{E}$ being $\mathcal{E}$-connected is a stronger condition to impose than being source connected. However in [2] we show that a Lie groupoid is source path/simply connected iff it is $\mathcal{E}$-path/$\mathcal{E}$-simply connected.

The proof of Lie's Second Theorem in Section 8 doesn't rely on the topological or smooth structure of the unit interval. In fact we can replace the pair groupoid on the unit interval with an arbitrary category $\mathbb{I}$ and a choice of arrow $\mathbf{2} \to \mathbb{I}$ in $Cat(\mathcal{E})$ where $\mathcal{E}$ is an arbitrary topos. It is easy to see that the $\mathcal{E}$-connectedness conditions can be reformulated using this data. In Section 7.3 and Section 7.4 we use these generalised connectedness conditions to describe how to express maps out of a simply connected category in terms of its 1- and 2-dimensional path spaces.

## 7.2 Truncated Cubical Objects

We arrange the (infinitesimal and macroscopic) path spaces into truncated cubical objects. Recall that the *2-truncated cube category* $\Box_2$ is the subcategory of $Man$ generated by the following arrows:

$$I^2 \;\substack{\leftarrow (1_I,0) \longrightarrow \\ \xrightarrow{\pi_1} \\ \leftarrow (1_I,1) \longrightarrow \\ \leftarrow (0,1_I) \longrightarrow \\ \xrightarrow{\pi_2} \\ \leftarrow (1,1_I) \longrightarrow}\; I \;\substack{\xleftarrow{1} \\ \xrightarrow{!} \\ \xleftarrow{0}}\; 1 \tag{4}$$

where $I$ is the unit interval. Recall that *the category $c_2\mathcal{E}$ of 2-truncated cubical objects in a category $\mathcal{E}$* is the functor category $[\Box_2^{op}, \mathcal{E}]$. The arrows of $c_2\mathcal{E}$ will be called *2-cubical maps*.

Let $\mathbb{I}$ be a category in $\mathcal{E}$ and $l : \mathbf{2} \to \mathbb{I}$ be an arbitrarily chosen arrow in $\mathbb{I}$. Then precomposing $l$ with the source and target inclusions $s, t : 1 \to \mathbf{2}$ gives arrows $0, 1 : 1 \to \mathbb{I}$. Hence

$$\mathbb{I}^2 \;\substack{\leftarrow (1_\mathbb{I},0) \longrightarrow \\ \xrightarrow{\pi_1} \\ \leftarrow (1_\mathbb{I},1) \longrightarrow \\ \leftarrow (0,1_\mathbb{I}) \longrightarrow \\ \xrightarrow{\pi_2} \\ \leftarrow (1,1_\mathbb{I}) \longrightarrow}\; \mathbb{I} \;\substack{\xleftarrow{1} \\ \xrightarrow{!} \\ \xleftarrow{0}}\; 1 \tag{5}$$

defines a functor $\Box_2 \to \mathcal{E}$. Let $\mathbb{C}$ be a category in $\mathcal{E}$. Then mapping into $\mathbb{C}$ determines a 2-truncated cubical object

$$\mathbb{C}^{\mathbb{I}^2} \rightrightarrows \mathbb{C}^{\mathbb{I}} \rightrightarrows \mathbb{C}^1$$

in $\mathcal{E}$ which we will call the *path 2-cubical object of $\mathcal{E}$*. Similarly

$$\mathbb{I}_\infty^2 \;\substack{\leftarrow (1_{\mathbb{I}_\infty},0) \longrightarrow \\ \xrightarrow{\pi_1} \\ \leftarrow (1_{\mathbb{I}_\infty},1) \longrightarrow \\ \leftarrow (0,1_{\mathbb{I}_\infty}) \longrightarrow \\ \xrightarrow{\pi_2} \\ \leftarrow (1,1_{\mathbb{I}_\infty}) \longrightarrow}\; \mathbb{I}_\infty \;\substack{\xleftarrow{1} \\ \xrightarrow{!} \\ \xleftarrow{0}}\; 1 \tag{6}$$

defines a functor $\Box_2 \to \mathcal{E}$. Then mapping into $\mathbb{C}$ determines a 2-truncated cubical object

$$\mathbb{C}^{\mathbb{I}_\infty^2} \rightrightarrows \mathbb{C}^{\mathbb{I}_\infty} \rightrightarrows \mathbb{C}^1$$

in $\mathcal{E}$ which we will call the *Weinstein 2-cubical object of $\mathcal{E}$*.

## 7.3 The Arrow Space of a Simply Connected Category

We now express the arrow space of a simply connected category $\mathbb{C}$ in terms of the paths and homotopies in $\mathbb{C}$.

**Notation 7.2.** We write $2_*\mathbb{I}$ for the pushout $\mathbb{I}\ {}_1+_0\ \mathbb{I}$. Similarly given an arrow $\Psi : \mathbb{C} \to \mathbb{X}$ we write $2_\times \Psi$ for the arrow $(\Psi\pi_1, \Psi\pi_2) : \mathbb{C}^{2_*\mathbb{I}} \to \mathbb{X}^{2_*\mathbb{I}}$. The category $\partial\mathbb{I}^2$ is the pushout

$$\begin{array}{ccc} \mathbf{2} & \xrightarrow{\delta} & 2_*\mathbb{I} \\ {\scriptstyle \delta}\downarrow & & \downarrow{\scriptstyle \iota_1} \\ 2_*\mathbb{I} & \xrightarrow{\iota_2} & \partial\mathbb{I}^2 \end{array}$$

in $Cat(\mathcal{E})$ where $\delta(l) = (\iota_1 l \circ_{2_*\mathbb{I}} \iota_2 l)$. We write $\iota : \partial\mathbb{I}^2 \to \mathbb{I}^2$ for the inclusion induced by the arrows $((0, 1_\mathbb{I}), (1_\mathbb{I}, 1))$ and $((1_\mathbb{I}, 0), (1, 1_\mathbb{I}))$. In this section the we will use the exponential notation (e.g. $\mathbb{C}^\mathbb{I}$) to denote the $\mathcal{E}$-valued hom-object.

**Lemma 7.3.** *For all simply connected categories $\mathbb{C}$ in $\mathcal{E}$ the diagram*

$$\mathbb{C}^{\mathbb{I}^2} \xrightarrow[\mathbb{C}^{\iota\iota_2}]{\mathbb{C}^{\iota\iota_1}} \mathbb{C}^{2_*\mathbb{I}} \xrightarrow{\mathbb{C}^\delta} \mathbb{C}^{\mathbf{2}}$$

*is a coequaliser in $\mathcal{E}$.*

*Proof.* Since $\mathcal{E}$ is a topos the arrow $\mathbb{C}^\delta$ is the coequaliser of its kernel pair. Hence

$$\mathbb{C}^{\partial\mathbb{I}^2} \xrightarrow[\mathbb{C}^{\iota_2}]{\mathbb{C}^{\iota_1}} \mathbb{C}^{2_*\mathbb{I}} \xrightarrow{\mathbb{C}^\delta} \mathbb{C}^{\mathbf{2}}$$

is a coequaliser in $\mathcal{E}$. But now the result follows from the hypothesis that $\mathbb{C}$ is simply connected. □

## 7.4 Mapping out of Simply Connected Categories

In this section we show that maps out of a simply connected category are completely determined by maps between their truncated path cubical objects. More precisely we will prove the following proposition:

**Proposition 7.4.** *Let $\mathbb{C}$ and $\mathbb{X}$ be categories where $\mathbb{C}$ is simply connected and*

$$\begin{array}{ccc} \mathbb{C}^{\mathbb{I}^2} \rightrightarrows \mathbb{C}^\mathbb{I} \rightrightarrows \mathbb{C}^1 \\ {\scriptstyle \Psi_2}\downarrow \quad {\scriptstyle \Psi_1}\downarrow \quad {\scriptstyle \Psi_0}\downarrow \\ \mathbb{X}^{\mathbb{I}^2} \rightrightarrows \mathbb{X}^\mathbb{I} \rightrightarrows \mathbb{X}^1 \end{array}$$

*is a 2-cubical map as defined in Section 7.2. Then there is a functor $\psi : \mathbb{C} \to \mathbb{X}$ with object map $\psi_0 = \Psi_0$ and arrow map $\psi_1$ satisfying $\psi_1\mathbb{C}^l = \mathbb{X}^l\Psi_1$.*

*Proof.* The functor $\psi : \mathbb{C} \to \mathbb{X}$ will have object map $\psi_0 = \Psi_0$ and arrow map $\psi_1$ given by the factorisation

$$\begin{array}{ccc} \mathbb{C}^{\mathbb{I}^2} \xrightarrow[\mathbb{C}^{\iota\iota_2}]{\mathbb{C}^{\iota\iota_1}} & \mathbb{C}^{2_*\mathbb{I}} \xrightarrow{\mathbb{C}^\delta} & \mathbb{C}^{\mathbf{2}} \\ {\scriptstyle \Psi_2}\downarrow & {\scriptstyle 2_\times\Psi_1}\downarrow & \downarrow{\scriptstyle \psi_1} \\ \mathbb{X}^{\mathbb{I}^2} \xrightarrow[\mathbb{X}^{\iota\iota_2}]{\mathbb{X}^{\iota\iota_1}} & \mathbb{X}^{2_*\mathbb{I}} \xrightarrow{\mathbb{X}^\delta} & \mathbb{X}^{\mathbf{2}} \end{array}$$

where the top line is a coequaliser by Lemma 7.3. The left square of

$$\begin{array}{ccc} \mathbb{C}^{\mathbb{I}} \xrightarrow{(1_{\mathbb{C}^{\mathbb{I}}}, \mathbb{C}^{!0})} \mathbb{C}^{2*\mathbb{I}} \xrightarrow{\mathbb{C}^{\delta}} \mathbb{C}^{\mathbf{2}} \\ \Psi_1 \downarrow \quad\quad 2_\times \Psi_1 \downarrow \quad\quad \psi_1 \downarrow \\ \mathbb{X}^{\mathbb{I}} \xrightarrow{(1_{\mathbb{X}^{\mathbb{I}}}, \mathbb{X}^{!0})} \mathbb{X}^{2*\mathbb{I}} \xrightarrow{\mathbb{X}^{\delta}} \mathbb{X}^{\mathbf{2}} \end{array} \qquad (7)$$

commutes because $\Psi_1$ is part of a 2-cubical map. The right square commutes by the definition of $\psi_1$. Hence $\psi_1 \mathbb{C}^l = \mathbb{X}^l \Psi_1$.

*Quotient Map is a Reflexive Graph Homomorphism.* Firstly the outer rectangle of

$$\begin{array}{ccc} \mathbb{C}^{\mathbb{I}} \xrightarrow{\mathbb{C}^l} \mathbb{C}^{\mathbf{2}} \xrightarrow[\mathbb{C}^0]{\mathbb{C}^1} \mathbb{C}^2 \\ \Psi_1 \downarrow \quad \psi_1 \downarrow \quad \Psi_0 \downarrow \\ \mathbb{X}^{\mathbb{I}} \xrightarrow{\mathbb{X}^l} \mathbb{X}^{\mathbf{2}} \xrightarrow[\mathbb{X}^0]{\mathbb{X}^1} \mathbb{X}^1 \end{array}$$

is serially commutative because $\Psi$ is a 2-cubical map. The left square commutes by (7). Therefore the right square is serially commutative because $\mathbb{C}^l$ is an epimorphism.

Secondly the right square of

$$\begin{array}{ccc} \mathbb{C}^1 \xrightarrow{\mathbb{C}^!} \mathbb{C}^{\mathbb{I}} \xrightarrow{\mathbb{C}^l} \mathbb{C}^{\mathbf{2}} \\ \Psi_0 \downarrow \quad \Psi_1 \downarrow \quad \psi_1 \downarrow \\ \mathbb{X}^1 \xrightarrow{\mathbb{X}^!} \mathbb{X}^{\mathbb{I}} \xrightarrow{\mathbb{X}^l} \mathbb{X}^{\mathbf{2}} \end{array}$$

commutes by (7). The left square commutes because $\Psi$ is a 2-cubical map. Therefore $\psi_1$ is a reflexive graph homomorphism.

*Quotient Map Preserves Composition.* It follows from (7) that the left square in

$$\begin{array}{ccc} \mathbb{C}^{2*\mathbb{I}} \xrightarrow{\mathbb{C}^{2*l}} \mathbb{C}^{2*\mathbf{2}} \xrightarrow{\mu_{\mathbb{C}}} \mathbb{C}^{\mathbf{2}} \\ 2_\times \Psi_1 \downarrow \quad 2_\times \psi_1 \downarrow \quad \psi_1 \downarrow \\ \mathbb{X}^{2*\mathbb{I}} \xrightarrow{\mathbb{X}^{2*l}} \mathbb{X}^{2*\mathbf{2}} \xrightarrow{\mu_{\mathbb{X}}} \mathbb{X}^{\mathbf{2}} \end{array}$$

commutes. Now $\mu_{\mathbb{C}} \mathbb{C}^{2*l} = \mathbb{C}^{\delta}$ and $\mu_{\mathbb{X}} \mathbb{X}^{2*l} = \mathbb{X}^{\delta}$ and so the outer square commutes by the definition of the quotient map $\psi_1$. Therefore the right hand square commutes because $\mathbb{C}^{2*l}$ is an epimorphism. Hence $\psi_1$ is an internal functor. □

## 7.5 Integrating Homomorphisms using Path Spaces

In this section we will show that in order to integrate homomorphisms out of a simply connected category $\mathbb{C}$ with path connected jet part it suffices to integrate the paths and homotopies in $\mathbb{C}$. More precisely, we will prove the following result:

**Proposition 7.5.** *Let $\mathbb{C}$ be a simply connected category in $\mathcal{E}$ such that the jet part $\mathbb{C}_\infty$ is path connected. Then any commutative square of the form*

$$\begin{array}{ccc} \mathbb{C}_\infty & \xrightarrow{\phi} & \mathbb{X} \\ {\scriptstyle \iota_\mathbb{C}^\infty}\downarrow & {\scriptstyle \psi}\nearrow & \downarrow{\scriptstyle r} \\ \mathbb{C} & \xrightarrow{\xi} & \mathbb{Y} \end{array} \qquad (8)$$

*has a filler $\psi$ iff for $n \in \{0, 1, 2\}$ the squares*

$$\begin{array}{ccc} \mathbb{C}_\infty^{\mathbb{I}^n} & \xrightarrow{\phi^{\mathbb{I}^n}} & \mathbb{X}^{\mathbb{I}^n} \\ {\scriptstyle (\iota_\mathbb{C}^\infty)^{\mathbb{I}^n}}\downarrow & {\scriptstyle \Psi_n}\nearrow & \downarrow{\scriptstyle r^{\mathbb{I}^n}} \\ \mathbb{C}^{\mathbb{I}^n} & \xrightarrow{\xi^{\mathbb{I}^n}} & \mathbb{Y}^{\mathbb{I}^n} \end{array}$$

*have fillers $\Psi_n$ that are components for a 2-cubical map.*

*Proof.* Suppose that $\Psi_n$ satisfies the above conditions. By Proposition 7.4 we obtain a functor $\psi : \mathbb{C} \to \mathbb{X}$ with object map $\psi_0 = \Psi_0$ and arrow map $\psi_1$ satisfying $\psi_1 \mathbb{C}^l = \mathbb{X}^l \Psi_1$. We now check that $\psi$ is a filler for (8). Firstly

$$r\psi_1 \mathbb{C}^l = r\mathbb{X}^l \Psi_1 = \mathbb{Y}^l r^{\mathbb{I}} \Psi_1 = \mathbb{Y}^l \xi^{\mathbb{I}} = \xi \mathbb{C}^l$$

so $r\psi_1 = \xi$ because $\mathbb{C}^l$ is an epimorphism. Secondly

$$\psi_1 (\iota_\mathbb{C}^\infty)^{\mathbf{2}} \mathbb{C}_\infty^l = \psi_1 \mathbb{C}^l (\iota_\mathbb{C}^\infty)^{\mathbb{I}} = \mathbb{X}^l \Psi_1 (\iota_\mathbb{C}^\infty)^{\mathbb{I}} = \mathbb{X}^l \phi^{\mathbb{I}} = \phi \mathbb{C}_\infty^l$$

so $\psi \iota_\mathbb{C}^\infty = \phi$ because $\mathbb{C}_\infty^l$ is an epimorphism. $\square$

# 8 Lie's Second Theorem

In this section we will formulate and prove Lie's second theorem relative to a category of local models. In Section 8.1 we prove the more general result that all of the jet part inclusions $\iota_\mathbb{C}^\infty : \mathbb{C}_\infty \rightarrowtail \mathbb{C}$ are in the left class of the integral factorisation system. Then in Section 8.2 we see how to deduce Lie's second theorem relative to a category of local models from this more general result.

## 8.1 Infinitesimal Inclusions are in Left Class

Now we will prove the fundamental lifting property involved in Lie's second theorem. More explicitly we will prove the following theorem.

**Theorem 8.1.** *Let $\mathbb{C}$ be a simply connected category in $\mathcal{E}$ such that the jet part $\mathbb{C}$ is path connected. Then $\iota_\mathbb{C}^\infty : \mathbb{C}_\infty \to \mathbb{C}$ is in the left class of the integral factorisation system. In other words, for all $r$ in the right class of the integral* *factorisation system and commutative diagrams*

$$\begin{array}{ccc} \mathbb{C}_\infty & \xrightarrow{\phi} & \mathbb{X} \\ {\iota_\mathbb{C}^\infty}\downarrow & \overset{\psi}{\nearrow} & \downarrow r \\ \mathbb{C} & \xrightarrow{\xi} & \mathbb{Y} \end{array}$$

*there is a unique filler $\psi$.*

*Proof. Existence of Solutions.* By Proposition 7.5 it will suffice to find for all $n \in \{0,1,2\}$ fillers $\Psi_n$ making

$$\begin{array}{ccc} \mathbb{C}_\infty^{\mathbb{I}^n} & \xrightarrow{\phi^{\mathbb{I}^n}} & \mathbb{X}^{\mathbb{I}^n} \\ {(\iota_\mathbb{C}^\infty)^{\mathbb{I}^n}}\downarrow & \overset{\Psi_n}{\nearrow} & \downarrow r^{\mathbb{I}^n} \\ \mathbb{C}^{\mathbb{I}^n} & \xrightarrow{\xi^{\mathbb{I}^n}} & \mathbb{Y}^{\mathcal{I}^n} \end{array} \qquad (9)$$

commute and which satisfy the relations defining a 2-cubical map. Now Remark 6.3 tells us that the right square in

$$(10)$$

is a pullback. Moreover by Lemma 4.8 the left arrow $(\iota_\mathbb{C}^\infty)^{\mathbb{I}_\infty^n}$ is invertible and so we can define $\Psi_n$ as the factorisation induced by the pair

$$\left(\xi^{\mathbb{I}^n}, \phi^{\mathbb{I}_\infty^n}\left((\iota_\mathbb{C}^\infty)^{\mathbb{I}_\infty^n}\right)^{-1}\mathbb{C}^{\iota_\mathbb{I}^\infty}\right)$$

We now check that (9) commutes. It is immediate that $r^{\mathbb{I}^n}\Psi_n = \xi^{\mathbb{I}^n}$. Finally we read off the equalities

$$\mathbb{X}^{(\iota_\mathbb{I}^\infty)^n}\Psi_n(\iota_\mathbb{C}^\infty)^{\mathbb{I}^n} = \mathbb{X}^{(\iota_\mathbb{I}^\infty)^n}\phi^{\mathbb{I}^n}$$

and

$$r^{\mathbb{I}^n}\Psi_n(\iota_\mathbb{C}^\infty)^{\mathbb{I}^n} = r^{\mathbb{I}^n}\phi^{\mathbb{I}^n}$$

from (10) and conclude that $\Psi_n(\iota_\mathbb{C}^\infty)^{\mathbb{I}^n} = \phi^{\mathbb{I}^n}$.

*Uniqueness of Solutions.* Let $\psi$ and $\chi$ be two functors making

$$\begin{array}{ccc} \mathbb{C}_\infty & \xrightarrow{\phi} & \mathbb{X} \\ \iota_{\mathbb{C}}^\infty \downarrow & \overset{\psi}{\underset{\chi}{\rightrightarrows}} & \downarrow r \\ \mathbb{C} & \xrightarrow{\xi} & \mathbb{Y} \end{array}$$

commute. We will show that $\psi = \chi$. First we note that it will suffice to show that $\psi^{\mathbb{I}} = \chi^{\mathbb{I}}$ because then

$$\psi^{\mathbf{2}} \mathbb{C}^l = \mathbb{X}^l \psi^{\mathbb{I}} = \mathbb{X}^l \chi^{\mathbb{I}} = \chi^{\mathbf{2}} \mathbb{C}^l$$

and $\mathbb{C}^l$ is an epimorphism. Furthermore since $\mathbb{X}^{\mathbb{I}^n}$ is a pullback it will suffice to check that

$$r^{\mathbb{I}} \psi^{\mathbb{I}} = \xi^{\mathbb{I}} = r^{\mathbb{I}} \chi^{\mathbb{I}}$$

and

$$\mathbb{X}^{\iota_{\mathbb{I}}^\infty} \psi^{\mathbb{I}} = \psi^{\mathbb{I}\infty} \mathbb{C}^{\iota_{\mathbb{I}}^\infty} = \phi^{\mathbb{I}\infty} (\iota_{\mathbb{C}}^\infty)^{\mathbb{I}\infty} \mathbb{C}^{\iota_{\mathbb{I}}^\infty} = \chi^{\mathbb{I}\infty} \mathbb{C}^{\iota_{\mathbb{I}}^\infty} = \mathbb{X}^{\iota_{\mathbb{I}}^\infty} \chi^{\mathbb{I}}$$

to conclude that $\psi = \chi$.

$\square$

## 8.2 Lie's Second Theorem

In this section we describe how our previous work allows us to prove Lie's second theorem relative to a category of local models. Recall that in Section 6.2 we constructed an adjunction

$$Cat_\infty(\mathcal{E}) \underset{(-)_\infty}{\overset{(-)_{int}}{\rightleftarrows}} \bot \; Cat_{int}(\mathcal{E})$$

which is analogous to the adjunction between the category of Lie groups and the category of formal group laws.

**Definition 8.2.** The category $Cat_{int}^{sc}(\mathcal{E})$ is the full subcategory of $Cat(\mathcal{E})$ whose objects are simply connected categories $\mathbb{C}$ such that $\mathbb{C}_\infty$ is path connected and the arrow

$$\mathbb{C}^{\mathbb{I}} \xrightarrow{\mathbb{C}^{\iota_{\mathbb{I}}^\infty}} \mathbb{C}^{\mathbb{I}_\infty}$$

is an isomorphism in $Cat(\mathcal{E})$.

**Remark 8.3.** The category $Cat_{int}^{sc}(\mathcal{E})$ is analogous to the category of simply connected Lie groups. Finally we record the result that is analogous to Lie's second theorem in this context.

**Corollary 8.4.** *The restriction of the functor* $(-)_\infty$ *to* $Cat_{int}^{sc}(\mathcal{E})$ *is full and faithful.*

*Proof.* Let $\mathbb{C}$ be a simply connected category in $\mathcal{E}$ such that $\mathbb{C}_\infty$ is path connected. Let $\mathbb{X}$ be a category in $\mathcal{E}$ such that $\mathbb{X}^{\iota_\mathbb{1}^\infty}$ is an isomorphism. This means that the arrow $! : \mathbb{X} \to 1$ is in the right class of the integral factorisation system. To see that $(-)_\infty$ is faithful let $\psi, \psi' : \mathbb{C} \to \mathbb{X}$ be internal functors such that $\psi_\infty = \psi'_\infty$. But then both $\psi$ and $\psi'$ are fillers for the square

$$\begin{array}{ccc} \mathbb{C}_\infty & \xrightarrow{\iota_\mathbb{X}^\infty \psi_\infty} & \mathbb{X} \\ \iota_\mathbb{C}^\infty \downarrow & & \downarrow ! \\ \mathbb{C} & \xrightarrow{!} & 1 \end{array}$$

and hence are equal by Theorem 8.1. To see that $(-)_\infty$ is full let $\phi : \mathbb{C}_\infty \to \mathbb{X}_\infty$ be an internal functor. Then by Theorem 8.1 the square

$$\begin{array}{ccc} \mathbb{C}_\infty & \xrightarrow{\iota_\mathbb{X}^\infty \phi} & \mathbb{X} \\ \iota_\mathbb{C}^\infty \downarrow & & \downarrow ! \\ \mathbb{C} & \xrightarrow{!} & 1 \end{array}$$

has a unique filler $\chi$ and $\chi_\infty = (\iota_\mathbb{C}^\infty)_\infty \chi_\infty = (\iota_\mathbb{X}^\infty)_\infty \phi_\infty = \phi_\infty$ as required. □

## Acknowledgements

The author would like to acknowledge the assistance of Richard Garner, my Ph.D. supervisor at Macquarie University Sydney, who provided valuable comments and insightful discussions in the genesis of this work. In addition the author is grateful for the support of an International Macquarie University Research Excellence Scholarship.

# References

[1] Steve Awodey. *Category theory*, volume 49 of *Oxford Logic Guides*. The Clarendon Press, Oxford University Press, New York, 2006.

[2] Matthew Burke. Ordinary connectedness implies enriched connectedness and integrability for lie groupoids. *Work in Progress*, 2016.

[3] C. Cassidy, M. Hébert, and G. M. Kelly. Reflective subcategories, localizations and factorization systems. *J. Austral. Math. Soc. Ser. A*, 38(3):287–329, 1985.

[4] Alberto S. Cattaneo, Benoit Dherin, and Giovanni Felder. Formal symplectic groupoid. *Comm. Math. Phys.*, 253(3):645–674, 2005.

[5] Marius Crainic and Rui Loja Fernandes. Integrability of Lie brackets. *Ann. of Math. (2)*, 157(2):575–620, 2003.

[6] Brian Day. On adjoint-functor factorisation. In *Category Seminar (Proc. Sem., Sydney, 1972/1973)*, pages 1–19. Lecture Notes in Math., Vol. 420. Springer, Berlin, 1974.

[7] Michiel Hazewinkel. *Formal groups and applications*, volume 78 of *Pure and Applied Mathematics*. Academic Press, Inc. [Harcourt Brace Jovanovich, Publishers], New York-London, 1978.

[8] P. T. Johnstone. *Sketches of an elephant: A topos theory compendium.*, volume 3.

[9] P. T. Johnstone. *Topos theory*. Academic Press [Harcourt Brace Jovanovich, Publishers], London-New York, 1977. London Mathematical Society Monographs, Vol. 10.

[10] G. M. Kelly. A unified treatment of transfinite constructions for free algebras, free monoids, colimits, associated sheaves, and so on. *Bull. Austral. Math. Soc.*, 22(1):1–83, 1980.

[11] G. M. Kelly. Basic concepts of enriched category theory. *Repr. Theory Appl. Categ.*, (10):vi+137, 2005. Reprint of the 1982 original [Cambridge Univ. Press, Cambridge; MR0651714].

[12] Anders Kock. *Synthetic differential geometry*, volume 333 of *London Mathematical Society Lecture Note Series*. Cambridge University Press, Cambridge, second edition, 2006.

[13] Rory Lucyshyn-Wright. Enriched factorisation systems. *Theory and Application of Categories*, 29(18):475–495, 2014.

[14] Saunders Mac Lane and Ieke Moerdijk. *Sheaves in geometry and logic*. Universitext. Springer-Verlag, New York, 1994. A first introduction to topos theory, Corrected reprint of the 1992 edition.

[15] Ieke Moerdijk and Gonzalo E. Reyes. *Models for smooth infinitesimal analysis*. Springer-Verlag, New York, 1991.

[16] Emily Riehl. *Categorical homotopy theory*, volume 24 of *New Mathematical Monographs*. Cambridge University Press, Cambridge, 2014.

[17] Jean-Pierre Serre. *Lie algebras and Lie groups*, volume 1500 of *Lecture Notes in Mathematics*. Springer-Verlag, Berlin, 2006. 1964 lectures given at Harvard University, Corrected fifth printing of the second (1992) edition.